\documentclass[12pt]{article}
\usepackage{amsmath,amsfonts,amssymb}
\usepackage[enableskew]{youngtab}
\usepackage{theorem}
\usepackage[usenames]{color}
\nonstopmode

\newtheorem{theorem}{Theorem}[section]
\newtheorem{proposition}[theorem]{Proposition}
\newtheorem{lemma}[theorem]{Lemma}
\newtheorem{definition}[theorem]{Definition}

\theorembodyfont{\rmfamily}
\newtheorem{example}[theorem]{Example}
\newtheorem{remark}[theorem]{Remark}

\begin{document}
$\,$\vspace{10mm}

\begin{center}
{\textsf {\huge Kirillov--Schilling--Shimozono bijection}}
\vspace{2mm}\\
{\textsf {\huge as energy functions of crystals}}
\vspace{15mm}\\
{\textsf {\Large Reiho Sakamoto}}
\vspace{2mm}\\
{\textsf {Department of Physics, Graduate School of Science,}}
\vspace{-1mm}\\
{\textsf {University of Tokyo, Hongo, Bunkyo-ku,}}
\vspace{-1mm}\\
{\textsf {Tokyo, 113-0033, Japan}}
\vspace{20mm}
\end{center}

\begin{abstract}
\noindent
The Kirillov--Schilling--Shimozono (KSS) bijection
appearing in theory of the Fermionic formula
gives an one to one correspondence between
the set of elements of tensor products
of the Kirillov--Reshetikhin crystals (called paths) and
the set of rigged configurations.
It is a generalization of Kerov--Kirillov--Reshetikhin
bijection and plays inverse scattering formalism
for the box-ball systems.
In this paper, we give an algebraic
reformulation of the KSS map from
the paths to rigged configurations,
using the combinatorial $R$ and energy functions of crystals.
It gives a characterization of the KSS
bijection as an intrinsic property of tensor products of crystals.
\end{abstract}
\pagebreak

\section{Introduction}
In proving the Fermionic formula,
the combinatorial bijection called
the Kirillov--Schilling--Shimozono (KSS) bijection
\cite{KSS,Sch,DS} plays the central role.
The KSS bijection is a natural generalization
of the original bijection
of Kerov--Kirillov--Reshetikhin \cite{KKR,KR},
and gives elaborated
combinatorial one to one correspondences
between the set of rigged configurations
and elements of tensor products
of the Kirillov--Reshetikhin crystals $B^{r,s}$
which we call {\it paths}:
$$\phi :\mathrm{path}\longmapsto
\mbox{rigged configuration}.$$
Here, we are focusing on the $A^{(1)}_n$ case,
and the Kirillov--Reshetikhin crystals are crystals
for irreducible finite dimensional modules
over quantum affine algebras \cite{Kas,KMN2}
indexed by Dynkin node $r\in\{1,2,\ldots,n\}$
and $s\in\mathbb{Z}_{>0}$.
Under the KSS bijection, the statistics called cocharge
is preserved and hence it gives proof of
the so called $X=M$ identity.
Here $M$ is called Fermionic formula,
and $X$ is the generating function of
energy functions over paths.
In the case we are considering,
$X$ contains the Kostka polynomial \cite{M}
as special case.
Reviews by Okado \cite{O} and Schilling \cite{Sch1}
give excellent introduction to the subject.
The reader can find an alternative approach for the
problem in Section 8 of the paper \cite{SW}.

Recently, a remarkable new feature of the KSS bijection
was discovered.
In the paper \cite{KOSTY}, the KSS bijection is
identified as the inverse scattering formalism
for the celebrated example of ultradiscrete integrable systems
called the box-ball systems \cite{TS,Tak,HHIKTT,FOY},
and as an outcome of the result,
algebraic alternative procedure for the
calculation of the map $\phi^{-1}$
from the rigged configurations
to the paths of the form $B^{1,s_1}\otimes\cdots\otimes
B^{1,s_L}$ is discovered \cite{Sak1}.
The latter formalism is used to derive explicit
piecewise linear formula for $\phi^{-1}$ in this case
as well as general solutions for the box-ball systems \cite{KSY}.
These formulas are written in terms of
the ultradiscrete limit of the tau functions.
Here, the tau functions are common apparatus in
algebraic structure theory of solitons (see, e.g.,
\cite{Hir,MJD} for textbook treatments)
and they satisfy the Hirota bilinear form.
Interestingly, although these tau functions come
from the KP hierarchy \cite{Sat,JM},
its expression contains the cocharge that appear in
the Fermionic formula.
This reveals unexpected link between the $X=M$ identities
and the soliton theory.
It is also interesting to note that if we restrict
these tau functions to the paths with periodicity
(the periodic box-ball systems \cite{YT}),
we obtain ultradiscrete limit of the classical
Riemann theta functions \cite{KS1,KS2}.

So it is natural to ask what the representation
theoretical origin of the KSS bijection is.
As seen previously,
for the case $\phi^{-1}$ the problem is settled in \cite{Sak1}
if the paths are elements of $B^{1,s_1}\otimes\cdots\otimes
B^{1,s_L}$, and the result has already brought many insights.
In this paper, we consider the case $\phi$
in full generality, i.e., we consider all elements of
$B^{r_1,s_1}\otimes\cdots\otimes B^{r_L,s_L}$.
As the result, the KSS map $\phi$ is reformulated
by algebraic language such as combinatorial $R$
and energy function of the crystal bases.
In particular, it seems that the result of \cite{Sak1}
do not admit straightforward generalizations, therefore
our result gives alternative approach to the problem
beyond the theory given in \cite{Sak1}.

Let us explain our result in more details.
Consider an element
$b\in B^{r_1,s_1}\otimes\cdots\otimes B^{r_L,s_L}$.
Assume that $b$ is an element of an affine crystal.
Then we consider the following isomorphism
(see the end of Section 2.1 for notations)
$$u_l[0]\otimes b_1[0]\otimes\cdots\otimes b_L[0]
\simeq b_1'[-e_1]\otimes\cdots\otimes b_L'[-e_L]\otimes
u_l^\ast [E_l]$$
under the isomorphism of affine combinatorial $R$.
Here $u_l\in B^{k,l}$ is the highest element
and $u_l^\ast\in B^{k,l}$.
We remark that this isomorphism is nothing but
time evolution of the box-ball systems \cite{HHIKTT,FOY}.
Then our main result (Theorem \ref{th:main}) states that
by taking differences of these $E_l$ or $e_i$'s,
we can derive alternative algorithm of the map $\phi$.
This is a generalization of the procedure given in \cite{Sak2}
where the case $A_1^{(1)}$ is treated.

The key to derive such a result is formula
connecting $E_l$ in the above
and shapes of Young diagrams of the
rigged configuration (Proposition \ref{prop:E=Q}).
Note that by virtue of the Yang--Baxter relation,
we have $E_l(b)=E_l(c)$ if $b$ and $c$ are isomorphic
under the classical part of the combinatorial $R$.
Therefore our scheme naturally explains
why Young diagrams of rigged configuration
do not change for two isomorphic paths $b\simeq c$.
This is a part of result \cite{KSS} (see Theorem \ref{th:KSS}),
which gives foundation for the inverse scattering formalism
in \cite{KOSTY}.
{}From more technical point of view, the original
combinatorial description of the map $\phi$
uses combinatorial notions such as the singular strings.
However, in our scheme, these combinatorial procedures
are precisely rephrased in computation of the energy function
Eq.(\ref{def:epsilon}).
Therefore we can say that our scheme simplifies
the rather mysterious concept of the KSS map $\phi$.
Since our scheme is purely algebraic, it will be
a very interesting future problem to give an algebraic alternative
reformulation of the theory given in \cite{KSS}.

The organization of this paper is as follows.
In Section 2, we recall basic facts
about the Kirillov--Reshetikhin crystals,
especially explicit description of the combinatorial
$R$ and energy function following Shimozono \cite{Shimo}.
In Section 3, we give our main result (Theorem \ref{th:main}).
In Section 4, we give a formula which translates
the KSS bijection into the energy function.
In Section 5, we prove Theorem \ref{th:main}.
In Appendix, we collect necessary facts from
the KSS bijection including its definition.

\section{Preliminaries}
\subsection{Kirillov--Reshetikhin crystal}
Let $W_s^{(r)}$ be a $U_q(\mathfrak{g})$ Kirillov--Reshetikhin
module, where we shall consider the case $\mathfrak{g}=A_n^{(1)}$.
The module $W_s^{(r)}$ is indexed by a Dynkin node
$r\in\{1,2,\ldots,n\}$ and $s\in\mathbb{Z}_{>0}$.
As a $U_q(A_n)$-module, $W_s^{(r)}$ is isomorphic to the irreducible
module corresponding to the partition $(s^r)$.
For arbitrary $r$ and $s$, the module $W_s^{(r)}$
is known to have a crystal base \cite{Kas,KMN2}, 
which we denote by $B^{r,s}$.
As a set, $B^{r,s}$ contains all the column strict
semi-standard Young tableaux
of depth $r$ and width $s$ over the alphabet $\{1,2,\ldots,n+1\}$.
Actions of the Kashiwara operators
$\tilde{e}_i$, $\tilde{f}_i$
for $i\neq 0$ coincide with the one described in \cite{KN},
and actions of $\tilde{e}_0$, $\tilde{f}_0$ are
determined by using the promotion operator \cite{Shimo}.
Since we do not use explicit forms of these operators,
we omit the details.
See \cite{O} for complements of this section.

For two crystals $B$ and $B'$,
one can define the tensor product
$B\otimes B'=\{b\otimes b'\mid b\in B,b'\in B'\}$.
On the tensor product,
actions of the Kashiwara operators
have simple form.
Namely, the operators 
$\tilde{e}_i,\tilde{f}_i$ act on $B\otimes B'$ by
\begin{eqnarray*}
\tilde{e}_i(b\otimes b')&=&\left\{
\begin{array}{ll}
\tilde{e}_i b\otimes b'&\mbox{ if }\varphi_i(b)\ge\varepsilon_i(b')\\
b\otimes \tilde{e}_i b'&\mbox{ if }\varphi_i(b) < \varepsilon_i(b'),
\end{array}\right. \\
\tilde{f}_i(b\otimes b')&=&\left\{
\begin{array}{ll}
\tilde{f}_i b\otimes b'&\mbox{ if }\varphi_i(b) > \varepsilon_i(b')\\
b\otimes \tilde{f}_i b'&\mbox{ if }\varphi_i(b)\le\varepsilon_i(b').
\end{array}\right. 
\end{eqnarray*}
Here we set 
$\varepsilon_i(b)=\max\{m\ge0\mid \tilde{e}_i^m b\ne0\}$
and $\varphi_i(b)=\max\{m\ge0\mid \tilde{f}_i^m b\ne0\}$.
We assume that $0\otimes b'$ and $b\otimes 0$ as $0$.
Then it is known that there is a unique crystal isomorphism
$R:B^{r,s}\otimes B^{r',s'}
\stackrel{\sim}{\rightarrow}B^{r',s'}\otimes B^{r,s}$.
We call this map (classical) combinatorial $R$
and usually write the map $R$ simply by $\simeq$.

Let us consider the affinization of the crystal $B$.
As the set, it is
\begin{equation}
\mathrm{Aff}(B)=\{b[d]\, |\, b\in B,\, d\in\mathbb{Z}\}.
\end{equation}
For the tensor product
$b[d]\otimes b'[d']\in
\mathrm{Aff}(B)\otimes\mathrm{Aff}(B')$,
we can lift the (classical)
combinatorial $R$ to affine case
as follows:
\begin{equation}\label{eq:affineR}
b[d]\otimes b'[d']\stackrel{R}{\simeq}
\tilde{b}'[d'-H(b\otimes b')]\otimes
\tilde{b}[d+H(b\otimes b')],
\end{equation}
where $b\otimes b'\simeq \tilde{b}'\otimes \tilde{b}$
is the isomorphism of (classical) combinatorial $R$.
The function $H(b\otimes b')$ is called the energy function.
We will give explicit forms of the combinatorial $R$ and
energy function in the subsequent sections.
Before closing this subsection,
we note the following important property.
\begin{proposition}\label{prop:YBeq}
The following Yang--Baxter equation holds on
$\mathrm{Aff}(B)\otimes\mathrm{Aff}(B')\otimes\mathrm{Aff}(B'')$:
\begin{equation}
(R\otimes1)(1\otimes R)(R\otimes1)=
(1\otimes R)(R\otimes1)(1\otimes R).
\end{equation}
\end{proposition}

\subsection{Schensted bumping algorithm}
In order to give the explicit form of the combinatorial $R$
and energy function, we recall the Schensted
bumping algorithm \cite{Schensted}.
Consider the semi-standard Young tableau $Y$
and positive integer $x$.
Then the insertion $(Y\leftarrow x)$
is obtained by the following steps:
\begin{enumerate}
\item
Insert $x$ in the first row of $Y$ either by
displacing the leftmost smallest number which is strictly
larger than $x$, or if no number is larger than $x$,
by adding $x$ on the right of the first row.

\item
If $x$ displaced a number from the first row,
then insert this number in the second row either by
displacing the leftmost smallest number which is strictly
larger than it or by adding it on the right of
the second row.

\item
Repeat this process row by row until
some number is added on the right of a row.
Here, we assume that there is an empty row
on the bottom of $Y$.
\end{enumerate}
We give an example of this procedure:
$$\left(
\Yvcentermath1
\young(1111223455,22333445,34556)
\longleftarrow 2\right)=
\Yvcentermath1
\newcommand{\niast}{2^\ast}
\newcommand{\sanast}{3^\ast}
\newcommand{\yonast}{4^\ast}
\newcommand{\goast}{5^\ast}
\young(111122\niast 455,22333\sanast 45,34\yonast 56,\goast).
$$
Here, the numbers with asterisks are the displaced ones.
We denote the successive applications of the row insertions
by $(Y\leftarrow xy)=((Y\leftarrow x)\leftarrow y)$ and so on.

We can easily infer the inverse of the bumping procedure
$(Y\leftarrow x)$.
Namely, starting from the node
$(Y\leftarrow x)\setminus Y (=:y)$,
we do inverse of row insertions until $x$
is ejected from $Y$.
Elementary steps in the inverse procedure are as follows:
start from integer $y$ and insert it to
row $w$ just above $y$.
Denote by $y'$ the rightmost integer of $w$ that is
strictly smaller than $y$, then new row $w'$ is the row
which is obtained by displacing $y'$ by $y$.
Insert $y'$ to the upper row and repeat
until the top row.

\subsection{Algorithm for combinatorial R and energy function}
We give the explicit description of the combinatorial $R$
and the energy function on $B^{r,s}\otimes B^{r',s'}$.
We begin with a few terminologies on Young tableaux.
Denote rows of a Young tableau $Y$ by $y_1,y_2,\ldots ,y_r$
from top to bottom.
Then row word $row(Y)$ is defined by concatenating rows as
$row(Y)=y_ry_{r-1}\ldots y_1$.
Let $x=(x_1,x_2,\ldots )$ and $y=(y_1,y_2,\ldots )$ be two partitions.
Then concatenation of $x$ and $y$ is the partition
$(x_1+y_1,x_2+y_2,\ldots )$.

\begin{proposition}[\cite{Shimo}]\label{prop:shimozono}
$b\otimes b'\in B^{r,s}\otimes B^{r',s'}$ is mapped to
$\tilde{b}'\otimes \tilde{b}\in B^{r',s'}\otimes B^{r,s}$
under the combinatorial $R$, i.e.,
\begin{equation}
b\otimes b'\stackrel{R}{\simeq}\tilde{b}'\otimes\tilde{b},
\end{equation}
if and only if
\begin{equation}
(b'\leftarrow row(b))=(\tilde{b}\leftarrow row(\tilde{b}')).
\end{equation}
Moreover, the energy function $H(b\otimes b')$ is given by
the number of nodes of $(b'\leftarrow row(b))$
outside the concatenation of partitions
$(s^r)$ and $({{s'}^{r'}})$.
\end{proposition}

In order to describe the
algorithm for finding $\tilde{b}$ and $\tilde{b}'$ from
the data $(b'\leftarrow row(b))$,
we introduce a terminology.
Let $Y$ be a tableau, and $Y'$ be a subset of $Y$
such that $Y'$ is also a tableau.
Consider the set theoretic subtraction $\theta =Y\setminus Y'$.
If the number of nodes contained in $\theta$ is $r$,
and if the number of nodes of $\theta$ contained in
each row is always 0 or 1,
then $\theta$ is called vertical $r$-strip.

Given a tableau
$Y=(b'\leftarrow row(b))$, let $Y'$ be the upper left
part of $Y$ whose shape is $(s^r)$.
We assign numbers from 1 to $r's'$
for each node contained in $\theta =Y\setminus Y'$
by the following procedure.
Let $\theta_1$ be the vertical $r'$-strip of $\theta$
as upper as possible.
For each node in $\theta_1$,
we assign numbers 1 through $r'$ from bottom to top.
Next we consider $\theta\setminus\theta_1$,
and find the vertical $r'$ strip $\theta_2$
by the same way.
Continue this procedure until all nodes of $\theta$
are assigned numbers up to $r's'$.
Then we apply inverse bumping procedure according
to the labeling of nodes in $\theta$.
Denote by $u_1$ the integer which is ejected
when we apply inverse bumping procedure starting
from the node with label 1.
Denote by $Y_1$ the tableau such that
$(Y_1\leftarrow u_1)=Y$.
Next we apply inverse bumping procedure
starting from the node of $Y_1$ labeled by 2,
and obtain the integer $u_2$ and tableau $Y_2$.
We do this procedure until we obtain $u_{r's'}$
and $Y_{r's'}$.
Finally, we have
\begin{equation}
\tilde{b}'=(\emptyset\leftarrow u_{r's'}u_{r's'-1}\cdots
u_1),\qquad
\tilde{b}=Y_{r's'}.
\end{equation}
\begin{example}
Consider the following tensor product:
$$b\otimes b'=
\Yvcentermath1
\young(11,24)\otimes\young(34,45,56)\in
B^{2,2}\otimes B^{3,2}.$$
{}From $b$, we have $row(b)=2411$, hence we have
$$\left(
\Yvcentermath1\young(34,45,56)
\leftarrow 2411
\right)=
\newcommand{\yonsan}{4_3}
\newcommand{\sanroku}{3_6}
\newcommand{\goni}{5_2}
\newcommand{\yongo}{4_5}
\newcommand{\rokuichi}{6_1}
\newcommand{\goyon}{5_4}
\Yvcentermath1
\young(11\yonsan ,24,\sanroku\goni ,\yongo\rokuichi ,\goyon) .
$$
Here subscripts of each node indicate the order of
inverse bumping procedure.
For example, we start from the node $6_1$ and obtain
$$\left(
\Yvcentermath1
\young(144,25,36,4,5)
\leftarrow 1
\right)=
\Yvcentermath1
\young(114,24,35,46,5),
\qquad\mathrm{therefore,}\qquad
Y_1=
\newcommand{\yonsan}{4_3}
\newcommand{\sanroku}{3_6}
\newcommand{\goni}{5_2}
\newcommand{\yongo}{4_5}
\newcommand{\rokuni}{6_2}
\newcommand{\goyon}{5_4}
\Yvcentermath1
\young(14\yonsan ,25,\sanroku\rokuni ,\yongo ,\goyon)
,\qquad
u_1=1.
$$
Next we start from the node $6_2$ of $Y_1$.
Continuing in this way, we obtain
$u_6u_5\cdots u_1=321541$ and
$Y_6=\Yvcentermath1
\young(44,56)$.
Since $(\emptyset\leftarrow 321541)=
\Yvcentermath1\young(11,24,35)$,
we obtain
$$\Yvcentermath1
\young(11,24)\otimes\young(34,45,56)\simeq
\young(11,24,35)\otimes\young(44,56)\, ,
\qquad
H\left(
\Yvcentermath1
\young(11,24)\otimes\young(34,45,56)
\right)
=3.
$$
Note that the energy function is derived from the
concatenation of shapes of $b$ and $b'$,
i.e., $\Yvcentermath1\yng(4,4,2)\,$.
\end{example}

\section{Kirillov--Schilling--Shimozono bijection
as energy function}
In this section, we present a procedure
to obtain the image of the Kirillov--Schilling--Shimozono
(KSS) bijection by using the energy function and combinatorial $R$.
Necessary facts about the KSS bijection are summarized in Appendix.
Let us consider the element $b=b_1\otimes b_2\otimes\cdots\otimes b_L
\in B^{\alpha_1,\beta_1}\otimes B^{\alpha_2,\beta_2}
\otimes\cdots \otimes B^{\alpha_L,\beta_L}$,
which we call path.
Considering $b_j$ as tableau, we denote each column of it
by $b_j=c_{\beta_j}\cdots c_2c_1$.
Then we define the subsets of the tableau
$b_{j,k}=c_k\cdots c_2c_1$.
We express the isomorphism of the combinatorial
$R$
\begin{equation}
a\otimes b\simeq b'\otimes a'
\end{equation}
by the following vertex diagram:
\begin{center}
\unitlength 13pt
\begin{picture}(4,4)
\put(0,2.0){\line(1,0){3.2}}
\put(1.6,1.0){\line(0,1){2}}
\put(-0.6,1.8){$a$}
\put(1.4,0){$b'$}
\put(1.4,3.2){$b$}
\put(3.4,1.8){$a'$}
\put(4.1,1.7){.}
\end{picture}
\end{center}
Successive applications of the combinatorial $R$
are depicted by concatenating these vertices.
Let $u_{l,0}^{(a)}=u_{l}^{(a)}\in B^{a,l}$ be the highest weight
element and define $u_{l,j}^{(a)}$ by the following diagram.
\begin{equation}\label{def:T_l}
\unitlength 13pt
\begin{picture}(22,5)(0,-0.5)
\multiput(0,0)(5.8,0){2}{
\put(0,2.0){\line(1,0){4}}
\put(2,0){\line(0,1){4}}
}
\put(-1.4,1.8){$u_{l,0}^{(a)}$}
\put(1.7,4.2){$b_1$}
\put(1.7,-0.8){$b_1'$}
\put(4.2,1.8){$u_{l,1}^{(a)}$}
\put(7.5,4.2){$b_2$}
\put(7.5,-0.8){$b_2'$}
\put(10.0,1.8){$u_{l,2}^{(a)}$}
\multiput(11.5,1.8)(0.3,0){10}{$\cdot$}
\put(14.7,1.8){$u_{l,L-1}^{(a)}$}
\put(17,0){
\put(0,2.0){\line(1,0){4}}
\put(2,0){\line(0,1){4}}
}
\put(18.7,4.2){$b_L$}
\put(18.7,-0.8){$b_L'$}
\put(21.2,1.8){$u_{l,L}^{(a)}$}
\end{picture}
\end{equation}
We call $u_{l,j}^{(a)}$ {\it carrier} following
terminology used in the box-ball systems.
Here the highest weight element
$u_{l}^{(a)}\in B^{a,l}$ is the tableau
whose $i$-th row is occupied by integers $i$.
For example, $u_4^{(3)}=\Yvcentermath1
\young(1111,2222,3333)\,$.
In the above diagram, let us write
$b'=b_1'\otimes b_2'\otimes\cdots\otimes b_L'$.
Then we define operator $T_l^{(a)}$ by $T_l^{(a)}(b)=b'$,
which is called
{\it time evolution operator of the box-ball systems.}

Associated with each vertex in (\ref{def:T_l}), we define
\begin{equation}\label{def:E_ljk}
E_{l,j,k}^{(a)}=H\left( u_{l,j-1}^{(a)}\otimes b_{j,k}\right),
\qquad (l\in\mathbb{Z}_{>0},
1\leq j\leq L, 1\leq k\leq \beta_j).
\end{equation}
We set $E_{0,j,k}^{(a)}=0$, $E_{l,j,0}^{(a)}=0$
and $E_{l,j,k}^{(0)}=0$.
\begin{definition}\label{def:led}
Local energy distribution (or energy spectrum)
is the collection of tables
such that the $l$-th row and the $(j,k)$-th column of
the $a$-th table
$(1\leq a\leq n$ for $A_n^{(1)})$ is given by
\begin{equation}\label{def:epsilon}
\varepsilon_{l,j,k}^{(a)}:=
\left( E_{l,j,k}^{(a)}-E_{l,j,k-1}^{(a)}\right) -
\left( E_{l-1,j,k}^{(a)}-E_{l-1,j,k-1}^{(a)}\right) .
\end{equation}
Here we assume that the numbering of the columns $(j,k)$
is according to the lexicographic ordering, i.e.,
$(p,q)<(j,k)$ if $p<j$ and $(j,q)\leq (j,k)$ if $q\leq k$.
\end{definition}
\begin{remark}\label{rem:meaning}
The meaning of column coordinates
$(j,k)$ in the above definition is that
$j$ specify the corresponding tensor factor $b_j$
and $k$ specify each column of tableau $b_j$.
We will explain this in more detail by the following diagram.
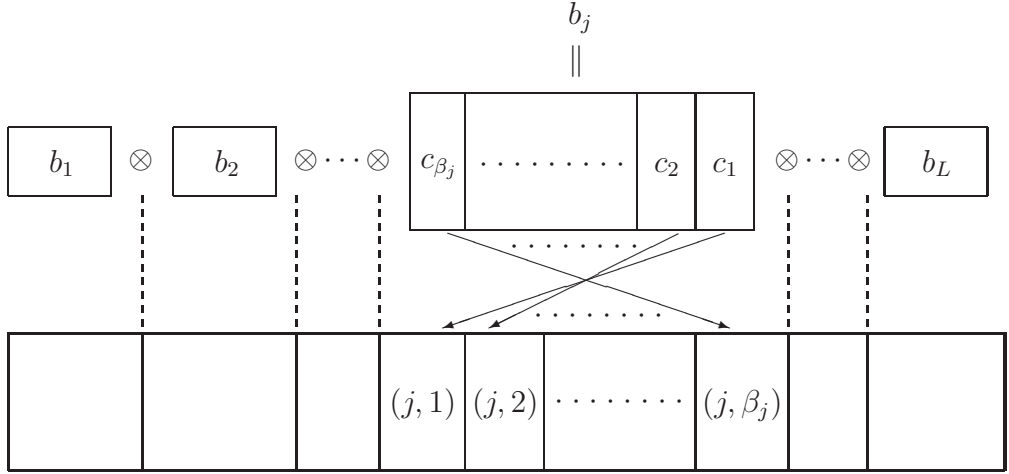
\begin{center}
\unitlength 13pt
\begin{picture}(29,13.5)
\multiput(0,8)(3,0){2}{\line(0,1){2}}
\multiput(0,8)(0,2){2}{\line(1,0){3}}
\put(1.2,8.7){$b_1$}
\put(3.5,8.8){$\otimes$}
\multiput(4.8,8)(3,0){2}{\line(0,1){2}}
\multiput(4.8,8)(0,2){2}{\line(1,0){3}}
\put(6.0,8.7){$b_2$}
\put(8.3,8.8){$\otimes\cdots\otimes$}
\multiput(11.7,7)(10,0){2}{\line(0,1){4}}
\multiput(11.7,7)(0,4){2}{\line(1,0){10}}
\put(16.3,11.7){$\|$}
\put(16.3,13){$b_j$}
\put(12,8.8){$c_{\beta_j}$}
\put(13.3,7){\line(0,1){4}}
\multiput(13.7,8.7)(0.5,0){9}{$\cdot$}
\multiput(18.3,7)(1.7,0){2}{\line(0,1){4}}
\put(20.5,8.7){$c_1$}
\put(18.8,8.7){$c_2$}
\put(22.3,8.8){$\otimes\cdots\otimes$}
\multiput(25.5,8)(3,0){2}{\line(0,1){2}}
\multiput(25.5,8)(0,2){2}{\line(1,0){3}}
\put(26.6,8.7){$b_L$}
\thicklines
\multiput(0,0)(29,0){2}{\line(0,1){4}}
\multiput(0,0)(0,4){2}{\line(1,0){29}}
\thinlines
\put(20.85,6.9){\vector(-3,-1){8.2}}
\put(19.5,6.9){\vector(-2,-1){5.5}}
\put(12.8,6.9){\vector(3,-1){8.2}}
\multiput(14.6,6.3)(0.5,0){8}{$\cdot$}
\multiput(15.3,4.3)(0.5,0){8}{$\cdot$}
\thicklines
\put(3.9,0){\line(0,1){4}}
\multiput(3.9,4.2)(0,0.4){10}{\line(0,1){0.2}}
\put(8.4,0){\line(0,1){4}}
\multiput(8.4,4.2)(0,0.4){10}{\line(0,1){0.2}}
\put(10.8,0){\line(0,1){4}}
\multiput(10.8,4.2)(0,0.4){10}{\line(0,1){0.2}}
\put(11.1,1.7){$(j,1)$}
\thinlines
\put(13.3,0){\line(0,1){4}}
\put(13.5,1.7){$(j,2)$}
\put(15.6,0){\line(0,1){4}}
\put(20.0,0){\line(0,1){4}}
\put(20.2,1.7){$(j,\beta_j)$}
\thicklines
\put(22.7,0){\line(0,1){4}}
\multiput(22.7,4.2)(0,0.4){10}{\line(0,1){0.2}}
\put(25.0,0){\line(0,1){4}}
\multiput(25.0,4.2)(0,0.4){10}{\line(0,1){0.2}}
\thinlines
\multiput(15.9,1.8)(0.5,0){8}{$\cdot$}
\end{picture}
\end{center}
In the above diagram, the first line shows the path
$b=b_1\otimes b_2\otimes\cdots\otimes b_L$
and the second line shows the corresponding
table for local energy distribution.
Columns of the local energy distributions are
labeled by $(1,1)$, $\cdots$, $(1,\beta_1)$,
$\cdots$, $(L,1)$, $\cdots$, $(L,\beta_L)$
from left to right.
In other words, the local energy distribution
is divided into $L$ regions corresponding to
$b_1$, $b_2$, $\cdots$, $b_L$ from left to right.
Now let us look at $b_j$ in more details.
Recall that we named columns of tableau $b_j$
by $b_j=c_{\beta_j}\cdots c_2c_1$, and
defined $b_{j,k}=c_k\cdots c_2c_1$.
As we see in Eq.(\ref{def:epsilon}),
$\varepsilon^{(a)}_{l,j,k}$ is defined by
taking difference between $k$ and $k-1$.
In other words, we consider difference between
two tableaux $b_{j,k}$ and $b_{j,k-1}$.
In this meaning, $(j,k)$-th column corresponds
to the column $c_k$.
Note that $c_k$ is the $k$-th column of $b_j$
from right,
on the other hand, $(j,k)$ is the $k$-th column
of the local energy distribution
(to be more specific, sub-region of it
corresponding to $b_j$) from left.
In this sense, left and right are reversed within
the region of the local energy distribution
corresponding to $b_j$.
\end{remark}

Given $b$, consider the following procedure
from Step 1 to 5.
\begin{enumerate}
\item
Draw the local energy distribution for $b$.
We do the following steps independently with respect to
the $a$-th table
($1\leq a\leq n$).

\item
Starting from the rightmost strictly
positive integer of the top row,
choose strictly positive integers successively as follows.
Assume that we have already chosen the one in the $l$-th row,
$k$-th column.
Take the set of all strictly positive integers of the $(l+1)$-th row
that are located
between the $k$-th column and the right end of the table.

As an explanation, we give a schematic diagram.
In the following diagram, rectangular frame represent
the local energy distribution and we
denote the chosen letter at $i$-th row by $e_i$.
Especially, the letter $e_l$ is at $l$-th row,
$k$-th column.
All other letters are suppressed for simplicity.
Then we consider the set of all strictly positive
integers contained in the gray strip at $(l+1)$-th row.
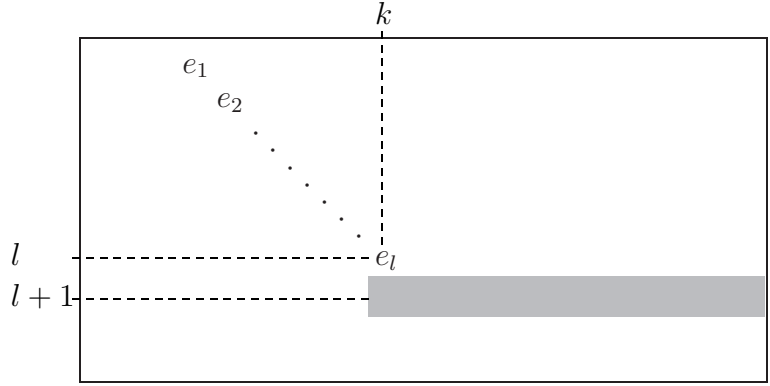
\begin{center}
\unitlength 13pt
\begin{picture}(20,10.5)
\multiput(0,0)(0,10){2}{\line(1,0){20}}
\multiput(0,0)(20,0){2}{\line(0,1){10}}
\put(3,9){$e_1$}
\put(4,8){$e_2$}
\multiput(5,7)(0.5,-0.5){7}{$\cdot$}
\put(8.6,3.4){$e_l$}
\color[cmyk]{0,0,0,0.3}
\put(8.4,1.9){\rule{150pt}{15pt}}
\color{black}
\multiput(8.8,4)(0,0.4){16}{\line(0,1){0.2}}
\put(8.6,10.4){$k$}
\multiput(-0.2,3.6)(0.4,0){22}{\line(1,0){0.2}}
\put(-2,3.4){$l$}
\multiput(-0.2,2.4)(0.4,0){22}{\line(1,0){0.2}}
\put(-2,2.2){$l+1$}
\end{picture}
\end{center}
If the set is empty, then stop.
Otherwise choose the rightmost element of the set and continue
until stop.
This procedure enables us to make the rightmost
possible successive group consisting of strictly positive integers.
For example, in the above diagram, letters
on the right of $e_i$ of the $i$-th row are all 0.
Write the position of the lastly picked integer
as $\mu_1^{(a)}$-th row, $(j_1^{(a)},k_1^{(a)})$-th column.

\item
Continue the previous grouping procedure
from right to left.
To be more precise, after making a group, we subtract
$1$ from all letters in the group and get new
table.
Then we apply Step 2 to this new table again
to get the next group.
We continue this grouping procedure until
all elements of the $a$-th table become $0$.
{}From lower endpoints of thus obtained groups,
we obtain the following sequence:
\begin{equation}
\left(\mu_1^{(a)},(j_1^{(a)},k_1^{(a)})\right),
\left(\mu_2^{(a)},(j_2^{(a)},k_2^{(a)})\right),
\cdots,
\left(\mu_{N^{(a)}}^{(a)},(j_{N^{(a)}}^{(a)},
k_{N^{(a)}}^{(a)})\right).
\end{equation}

\item
Calculate integers $r_1^{(a)}$, $r_2^{(a)}$, $\cdots$,
$r_{N^{(a)}}^{(a)}$ as follows:
\begin{eqnarray}
r_s^{(a)}&=&\mathcal{C}+\mathcal{E}
\label{eq:rigging1}\\
\mathcal{C}&=&\sum_{i=1}^{j_s^{(a)}-1}\delta_{\alpha_i,a}
\min \left( \mu_s^{(a)},\beta_i\right)
+\delta_{\alpha_{j_s^{(a)}},a}\min (\mu_s^{(a)},k_s^{(a)})
\label{eq:rigging2}\\
\mathcal{E}&=&\sum_{(p,q)\leq (j_s^{(a)},k_s^{(a)})}\,\,\,
\sum_{l\leq\mu_s^{(a)}}
\left(
\varepsilon_{l, p, q}^{(a-1)}-
2\varepsilon_{l, p, q}^{(a)}+
\varepsilon_{l, p, q}^{(a+1)}
\right)\label{eq:rigging3}
\end{eqnarray}
In other words,
summation of $\mathcal{E}$ runs over all entries of
upper left rectangular region of the local energy distribution
whose lower right corner is $\mu_s^{(a)}$-th row,
$(j_s^{(a)},k_s^{(a)})$-th column.

\item
{}From the shape of $b$, we extract the data
\begin{equation}
\mathcal{Q}=\left(
(\nu_i^{(0)})_{i=1}^{L^{(0)}},
(\nu_i^{(1)})_{i=1}^{L^{(1)}}, \cdots,
(\nu_i^{(n-1)})_{i=1}^{L^{(n-1)}}\right)
\end{equation}
as follows.
We set $\mathcal{Q}=\emptyset$ as the initial condition.
Then, proceeding from the left end of $b$ to right,
we add $\nu_i^{(a)}$ to $\mathcal{Q}$
if there is a factor of shape $B^{a+1,\nu_i^{(a)}}$.
Finally, we obtain the following data:
\begin{equation}
\left(
(\nu_i^{(0)})_{i=1}^{L^{(0)}},\cdots,
(\nu_i^{(n-1)})_{i=1}^{L^{(n-1)}},
(\mu_i^{(1)},r_i^{(1)})_{i=1}^{N^{(1)}},\cdots,
(\mu_i^{(n)},r_i^{(n)})_{i=1}^{N^{(n)}}
\right).
\end{equation}
\end{enumerate}
The following is the main result of this paper.
\begin{theorem}\label{th:main}
Given the path
$b=b_1\otimes b_2\otimes\cdots\otimes b_L
\in B^{\alpha_1,\beta_1}\otimes B^{\alpha_2,\beta_2}
\otimes\cdots \otimes B^{\alpha_L,\beta_L}$,
we obtain integers $\nu_i^{(a)}$, $(\mu_i^{(a)},r_i^{(a)})$
according to the above procedure from Step 1 to 5.
Then the result coincides with the image of the KSS bijection $\phi$:
\begin{equation}
\phi (b)=
\left(
(\nu_i^{(0)})_{i=1}^{L^{(0)}},\cdots,
(\nu_i^{(n-1)})_{i=1}^{L^{(n-1)}},
(\mu_i^{(1)},r_i^{(1)})_{i=1}^{N^{(1)}},\cdots,
(\mu_i^{(n)},r_i^{(n)})_{i=1}^{N^{(n)}}
\right).
\end{equation}
\end{theorem}
\begin{remark}
(a) If we start from non-highest weight element $b$,
the above theorem is true if we replace
the ``rigged configuration" by the
``unrestricted rigged configuration"
(see the latter part of Section \ref{app:facts}
for explanation).
(b) In the local energy distribution, the entries are
given by 0 or strictly positive integers
(as we will see at the beginning of the proof),
and the number of non-zero entries is finite,
which will be evaluated in Proposition \ref{prop:E=Q}.
\end{remark}
\begin{remark}\label{rem:ist}
It is known that the time evolution operator
$T_l^{(a)}$ can be linearized by the KSS bijection
\cite{KOSTY}.
When $\phi (b)$ is given by
\begin{equation}
\left(
(\nu_i^{(0)})_{i=1}^{L^{(0)}},\cdots,
(\nu_i^{(n-1)})_{i=1}^{L^{(n-1)}},
(\mu_i^{(1)},r_i^{(1)})_{i=1}^{N^{(1)}},\cdots,
(\mu_i^{(n)},r_i^{(n)})_{i=1}^{N^{(n)}}
\right),
\end{equation}
then $\phi \left(T_l^{(a)}(b)\right)$ is given by
\begin{align}
&\left(
(\nu_i^{(0)})_{i=1}^{L^{(0)}},\cdots,
(\nu_i^{(n-1)})_{i=1}^{L^{(n-1)}},\right.
\nonumber\\
&\qquad\left.
(\mu_i^{(1)},r_i^{(1)})_{i=1}^{N^{(1)}},\cdots,
(\mu_i^{(a)},r_i^{(a)}+\min (l,\mu_i^{(a)})
)_{i=1}^{N^{(a)}},\cdots,
(\mu_i^{(n)},r_i^{(n)})_{i=1}^{N^{(n)}}
\right).
\end{align}
In this sense, the groups obtained in Theorem
\ref{th:main} represent solitons contained in
a path.
The proof of this property uses Theorem \ref{th:KSS}.
For the tensor products of the form
$B^{1,s_1}\otimes\cdots\otimes
B^{1,s_L}$, these solitons actually identified
with those arise in the KP hierarchy \cite{KSY}.
\end{remark}

\begin{example}\label{ex:led}
Let us consider the following path:
\begin{equation}
b=\Yvcentermath1\young(1111)\otimes
\young(12,23,34)\otimes\young(1124,2235)
\end{equation}
Note that we give calculation according
to the original combinatorial procedure in
Example \ref{ex:KSS}.
If we compare Example \ref{ex:KSS} with
the data given in the following,
the reader will find that if we remove
a box from $l$-th column of $\mu^{(a)}$
then there is 1 at $l$-th row of
the corresponding column of $a$-th table of
the local energy distribution
(see Remark \ref{rem:meaning} for the correspondence
of columns of the path and the local energy distribution).

\pagebreak

Now the local energy distribution
(Definition \ref{def:led}) for $b$
takes the following forms.
$$a=1:\qquad
\begin{array}{|cccc|cc|cccc|}
\hline
0&0&0&0&1&0&0&0&0&0\\
0&0&0&0&0&0&1&0&0&0\\
0&0&0&0&0&0&0&1&0&0\\
0&0&0&0&0&0&0&0&0&0\\
\hline
\end{array}$$\\*[-6mm]
$$a=2:\qquad
\begin{array}{|cccc|cc|cccc|}
\hline
0&0&0&0&1&0&1&0&0&0\\
0&0&0&0&0&0&1&0&0&0\\
0&0&0&0&0&0&0&1&0&0\\
0&0&0&0&0&0&0&0&0&0\\
\hline
\end{array}$$\\*[-6mm]
$$a=3:\qquad
\begin{array}{|cccc|cc|cccc|}
\hline
0&0&0&0&1&0&1&0&0&0\\
0&0&0&0&0&0&1&0&0&0\\
0&0&0&0&0&0&0&0&0&0\\
\hline
\end{array}$$\\*[-6mm]
$$a=4:\qquad
\begin{array}{|cccc|cc|cccc|}
\hline
0&0&0&0&0&0&1&0&0&0\\
0&0&0&0&0&0&0&0&0&0\\
\hline
\end{array}$$
We classify the strictly positive integers contained
in the above tables as follows:
$$a=1:\qquad
\begin{array}{|cccc|cc|cccc|}
\hline
\quad\,\,\,\,& & & &1\,\,&\,\,& & & &\,\,\,\,\\
 & & & & & &1& & & \\
 & & & & & & &1^\ast& & \\
 & & & & & & & & & \\
\hline
\end{array}$$\\*[-6mm]
$$a=2:\qquad
\begin{array}{|cccc|cc|cccc|}
\hline
\quad\,\,\,\,& & & &2^\ast&\,\,&1& & &\,\,\,\,\\
 & & & & & &1& & & \\
 & & & & & & &1^\ast& & \\
 & & & & & & & & & \\
\hline
\end{array}$$\\*[-6mm]
$$a=3:\qquad
\begin{array}{|cccc|cc|cccc|}
\hline
\quad\,\,\,\,& & & &2^\ast&\,\,&1& & &\,\,\,\,\\
 & & & & & &1^\ast&\,\,\, & & \\
 & & & & & & & & & \\
\hline
\end{array}$$\\*[-6mm]
$$a=4:\qquad
\begin{array}{|cccc|cc|cccc|}
\hline
\quad\,\,\,\,& & & &\quad &\,\,&1^\ast&\,\,\, & &\,\,\,\,\\
 & & & & & & & & & \\
\hline
\end{array}$$
Here, according to the grouping procedure in Steps 2 and 3
in Theorem \ref{th:main},
labellings of each group are given from
right to left, and letters with asterisks
show lower end points of each group.

In the $a=1$-st table, we recognize one group of
cardinality 3 (labeled by 1)
whose end point is at
3-rd row, $(3,2)$-th column.
According to Eq.(\ref{eq:rigging2}), we have
$$\mathcal{C}=\min (3,4)+0+0=3,$$
where the first term in the middle
shows the contribution from \young(1111)
(i.e., $(\alpha_1,\beta_1)=(1,4)$).
According to Eq.(\ref{eq:rigging3}), we have
(the contribution from
$\varepsilon_{l,j,k}^{(0)}=0)$
$$\mathcal{E}=\{0+0+0+0\}+\{(-2\cdot 1+1)+0\}
+\{(-2\cdot 1+2)+(-2\cdot 1+1)\}=-2,$$
hence we have $\mathcal{C}+\mathcal{E}=1$.
Here, the parentheses $\{\cdot\}$ represent tensor factors
from the left one to right,
and the integers inside represent
$\sum_{l=1}^\mu\varepsilon_{l,j,k}^{(a)}$
where $\mu$ is the cardinality of the corresponding
group.
Recall that $\sum_{l=1}^\mu\varepsilon_{l,j,k}^{(a)}$ is
sum of integers within the first $\mu$ rows
of $(j,k)$-th column of the local energy distribution.

In the $a=2$-nd table, we recognize two groups
of cardinalities 3 and 1 (labeled by 
1 and 2, respectively).
Let us start from the right group whose
end point is at 3-rd row, $(3,2)$-th column.
Then we have
\begin{align}
\mathcal{C}=&\, 0+0+\min (3,2)=2,
\nonumber\\
\mathcal{E}=&\,\{0+0+0+0\}+\{(1-2\cdot 1+1)+0\}
+\{(1-2\cdot 2+2)+(1-2\cdot 1+0)\}
\nonumber\\
=&\, -2,\nonumber
\end{align}
hence we have the rigging $\mathcal{C}+\mathcal{E}=0$.
The left group has end point at
1-st row, $(2,1)$-st column, and we have
\begin{align}
\mathcal{C}=&\, 0+0=0,
\nonumber\\
\mathcal{E}=&\,\{0+0+0+0\}+\{(1-2\cdot 1+1)\}
=0\nonumber
\end{align}
hence we have the rigging $\mathcal{C}+\mathcal{E}=0$.

In the $a=3$-rd table, we recognize two groups
of cardinalities 2 and 1 (labeled by 1 and 2,
respectively).
Let us start from the right group whose
end point is at 2-nd row, $(3,1)$-th column.
Then we have
\begin{align}
\mathcal{C}=&\, 0+\min (2,2)+0=2,
\nonumber\\
\mathcal{E}=&\,\{0+0+0+0\}+\{(1-2\cdot 1+0)+0\}
+\{(2-2\cdot 2+1)\}
= -1,\nonumber
\end{align}
hence we have the rigging $\mathcal{C}+\mathcal{E}=0$.
Next, we consider the left group
whose end point is at 1-st row, $(2,1)$-th column.
Then we have
\begin{align}
\mathcal{C}=&\, 0+\min (3,1)=1,
\nonumber\\
\mathcal{E}=&\,\{0+0+0+0\}+\{(1-2\cdot 1+0)\}
= -1,\nonumber
\end{align}
hence we have the rigging $\mathcal{C}+\mathcal{E}=0$.

Finally, in the $a=4$-th table, we recognize
one group whose end point is 1-st row,
$(3,1)$-th column.
Then we have
\begin{align}
\mathcal{C}=&0+0+0=0,
\nonumber\\
\mathcal{E}=&\,\{0+0+0+0\}+\{1+0\}
+\{1-2\cdot 1\}
= 0,\nonumber
\end{align}
hence we have the rigging $\mathcal{C}+\mathcal{E}=0$.

As the result, we obtain the following
rigged configuration:
\begin{center}
\unitlength 13pt
\begin{picture}(26,6)
\Yboxdim13pt
\put(1,4){\young(\hfill\hfill\hfill\hfill)}
\put(2.5,5.5){$\nu^{(0)}$}
\put(7,1){\young(\hfill\hfill\hfill)}
\put(10.2,1.1){1}
\put(8,2.5){$\mu^{(1)}$}
\put(7,4){\young(\hfill\hfill\hfill\hfill)}
\put(8.5,5.5){$\nu^{(1)}$}
\put(13,0){\young(\hfill\hfill\hfill,\hfill)}
\put(14.2,0.1){0}
\put(16.2,1.1){0}
\put(14,2.5){$\mu^{(2)}$}
\put(13,4){\young(\hfill\hfill)}
\put(13.5,5.5){$\nu^{(2)}$}
\put(19,0){\young(\hfill\hfill,\hfill)}
\put(20.2,0.1){0}
\put(21.2,1.1){0}
\put(19.5,2.5){$\mu^{(3)}$}
\put(24,1){\young(\hfill)}
\put(25.2,1.1){0}
\put(24,2.5){$\mu^{(4)}$}
\end{picture}
\end{center}
Here we put the riggings
on the right of the corresponding rows.
\end{example}
\begin{example}
When the path is not a highest weight element,
we can use the same procedure.
Consider the following path:
\begin{equation*}
b=\Yvcentermath1\young(122,233)\otimes\young(11)\otimes
\young(1122,2333)\otimes\young(111,222,334)\otimes
\young(11233,23444)\otimes\young(1112,2223,3334)
\end{equation*}

Then its local energy distribution takes the
following forms.
$$a=1:\qquad
\begin{array}{|ccc|cc|cccc|ccc|ccccc|cccc|}
\hline
1&0&0&0&0&1&0&0&0&0&0&0&1&0&0&0&0&0&0&0&0\\
0&1&0&0&0&0&1&0&0&0&0&0&0&0&0&0&0&1&0&0&0\\
0&0&0&0&0&0&0&0&0&0&0&0&0&1&0&0&0&0&0&0&0\\
0&0&0&0&0&0&0&0&0&0&0&0&0&0&1&0&0&0&0&0&0\\
0&0&0&0&0&0&0&0&0&0&0&0&0&0&0&0&0&0&0&0&0\\
\hline
\end{array}$$\\*[-6mm]
$$a=2:\qquad
\begin{array}{|ccc|cc|cccc|ccc|ccccc|cccc|}
\hline
1&0&0&0&0&1&0&0&0&0&0&0&2&0&0&0&0&0&0&0&0\\
0&1&0&0&0&0&1&0&0&0&0&0&0&1&0&0&0&1&0&0&0\\
0&0&0&0&0&0&0&1&0&0&0&0&0&0&1&0&0&0&0&0&0\\
0&0&0&0&0&0&0&0&0&0&0&0&0&1&0&1&0&0&0&0&0\\
0&0&0&0&0&0&0&0&0&0&0&0&0&0&0&0&0&0&0&0&0\\
\hline
\end{array}$$\\*[-6mm]
$$a=3:\qquad
\begin{array}{|ccc|cc|cccc|ccc|ccccc|cccc|}
\hline
0&0&0&0&0&0&0&0&0&1&0&0&1&0&0&0&0&0&0&0&0\\
0&0&0&0&0&0&0&0&0&0&0&0&0&1&0&0&0&1&0&0&0\\
0&0&0&0&0&0&0&0&0&0&0&0&0&0&1&0&0&0&0&0&0\\
0&0&0&0&0&0&0&0&0&0&0&0&0&0&0&0&0&0&0&0&0\\
\hline
\end{array}$$
We classify the strictly positive integers into groups
in the following way:
$$a=1:\qquad
\begin{array}{|ccc|cc|cccc|ccc|ccccc|cccc|}
\hline
3& & & & &2& &      & & &\,\,\,& &1\,\,\,& & & & & & & & \\
 &3^*& & & & &2&\,\,\,& & &      & &   & & & & &1^*\!& & & \\
 & & & & & & &      & & &      & &   &2& & & & & & & \\
 & & & & & & &      & & &      & &   & &2^*&\,& & & & & \\
 & &\,&\quad& & & &\quad& & & &\quad&   & & & & & & & &\quad\,\\
\hline
\end{array}$$\\*[-6mm]
$$a=2:\qquad
\begin{array}{|ccc|cc|cccc|ccc|ccccc|cccc|}
\hline
4& & & & &3& & & & &\,\,\,& &21\!& & & & & & & & \\
 &4^*& & & & &3& & & & & & &2& & & &1^*\!& & & \\
 & & & & & & &3& & & & & & &2& & & & & & \\
 & & & & & & & & & & & & &3^*\!\!& &2^*\!& & & & & \\
 & &\,&\quad& & & &\quad& & & &\quad&   & & & & & & & &\quad\,\\
\hline
\end{array}$$\\*[-6mm]
$$a=3:\qquad
\begin{array}{|ccc|cc|cccc|ccc|ccccc|cccc|}
\hline
\quad& & & & &\quad\,\,\,& & & &2& & &1\,\,\,& & & & & & & & \\
     & & & & &           & & & & & & &   &2& & & &1^*\!& & & \\
     & & & & &           & & & & & & &   & &2^*&\,& & & & & \\
 & &\,\,\,&\quad& & & &\,\,\,& & & &\quad&   & & & & & & & &\quad\,\\
\hline
\end{array}$$
{}From these data, we obtain the following
unrestricted rigged configuration:
\begin{center}
\unitlength 13pt
\begin{picture}(22,10)
\put(0,8){\young(\hfill\hfill)}
\put(0.5,9.5){$\nu^{(0)}$}
\put(4,1){\young(\hfill\hfill\hfill\hfill,%
\hfill\hfill,\hfill\hfill)}
\put(6.2,1.2){$-2$}
\put(6.2,2.2){$-2$}
\put(8.2,3.2){$-2$}
\put(5.5,4.5){$\mu^{(1)}$}
\put(4,6){\young(\hfill\hfill\hfill\hfill\hfill,%
\hfill\hfill\hfill\hfill,\hfill\hfill\hfill)}
\put(6,9.5){$\nu^{(1)}$}
\put(11,0){\young(\hfill\hfill\hfill\hfill,%
\hfill\hfill\hfill\hfill,\hfill\hfill,\hfill\hfill)}
\put(13.3,0.2){$0$}
\put(13.3,1.2){$0$}
\put(15.3,2.2){$0$}
\put(15.3,3.2){$0$}
\put(12.5,4.5){$\mu^{(2)}$}
\put(11,7){\young(\hfill\hfill\hfill\hfill,%
\hfill\hfill\hfill)}
\put(12.5,9.5){$\nu^{(2)}$}
\put(18,2){\young(\hfill\hfill\hfill,\hfill\hfill)}
\put(20.3,2.2){$3$}
\put(21.3,3.2){$4$}
\put(19,4.5){$\mu^{(3)}$}
\end{picture}
\end{center}
\end{example}

\section{Energy function and shape of configuration}
\subsection{Formula}
Consider a path
$b=b_1\otimes b_2\otimes\cdots\otimes b_L
\in B^{\alpha_1,\beta_1}\otimes B^{\alpha_2,\beta_2}
\otimes\cdots \otimes B^{\alpha_L,\beta_L}$.
The goal of this section is to show the following property.
\begin{proposition}\label{prop:E=Q}
We have
\begin{equation}\label{eq:E=Q}
E_l^{(a)}=Q_l^{(a)},
\end{equation}
where we write $E_l^{(a)}=\sum_{j=1}^L E_{l,j,\beta_j}^{(a)}$,
and $Q_l^{(a)}$ is defined in Eq.(\ref{def:Q}).
\end{proposition}
The meaning of the above quantity $E_l^{(a)}$ is as follows.
Let us assume that $b$ is an element of
a tensor product of affine crystals.
Then from Eq.(\ref{eq:affineR}), we have
\begin{equation}
u_l^{(a)}[0]\otimes b\simeq
T_l^{(a)}(b)\otimes u^\ast [E_l^{(a)}].
\end{equation}

$E_l^{(a)}$ is a generalization of conserved quantities of
the box-ball systems introduced by
Fukuda--Okado--Yamada \cite{FOY}.
Their result stems from the Yang--Baxter relation
(Proposition \ref{prop:YBeq}).

\subsection{Lemmas about energy functions}
We begin with  basic lemmas about
the energy function.
The following is the immediate consequence of
the definition in Proposition \ref{prop:shimozono}.
\begin{lemma}\label{lem_vacuum*vacuum}
Consider the highest elements $u_s^{(r)}\in B^{r,s}$
and $u_{s'}^{(r')}\in B^{r',s'}$ for arbitrary
$r,s,r',s'$.
Then we have
\begin{equation}
u_{s}^{(r)}\otimes u_{s'}^{(r')}\simeq
u_{s'}^{(r')}\otimes u_{s}^{(r)},\qquad
H\left(
u_{s}^{(r)}\otimes u_{s'}^{(r')}
\right)=0.
\end{equation}
\end{lemma}
Furthermore, we have the following property.
\begin{lemma}\label{lem:energy_conservation}
For any path $b$ and any $l,k\in\mathbb{Z}_{>0}$,
$1\leq r,a\leq n$ we have
\begin{equation}
E_l^{(r)}\left(u_k^{(a)}\otimes b\right)=
E_l^{(r)}\left(b\otimes u_k^{(a)}\right)=E_l^{(r)}(b).
\end{equation}
\end{lemma}
{\bf Proof.} {}From Lemma \ref{lem_vacuum*vacuum}, we have
the equality between the first and the third quantities.
Therefore it is enough to show
\begin{eqnarray}\label{eq:v*u=0}
H\left(v\otimes u_k^{(a)}\right)=0
\end{eqnarray}
for arbitrary $v\in B^{r,s}$.
Denote the rows of $v$ as $v_1,\ldots,v_r$ from
top to bottom.
Since $v$ is semi-standard, the letters of $v_i$
are greater than or equal to $i$.
Consider the insertion $(u_k^{(a)}\leftarrow v_r)$.
Since letters of $v_r$ are greater than or equal to 1,
and since $v_r=v_{r,1}v_{r,2}\cdots v_{r,s}$ is a
non-decreasing integer sequence, we have
\begin{equation}
\left(u_k^{(a)}\leftarrow v_r\right)=
\Yvcentermath1
\Yboxdim18pt
\newcommand{\vrichi}{v_{r,1}}
\newcommand{\vrni}{v_{r,2}}
\newcommand{\vrs}{v_{r,s}}
\newcommand{\cdotsr}{\begin{picture}(2,2)
\multiput(-4,4)(4.7,0){3}{\circle*{1}}
\end{picture}}
\young(11\cdotsr 1\vrichi\vrni\cdotsr\vrs,22\cdotsr 2,%
\cdotsr\cdotsr\cdotsr\cdotsr,aa\cdotsr a)
\end{equation}
Assuming that $r>1$,
we insert the row $v_{r-1}$ into this
$(u_k^{(a)}\leftarrow v_r)$.
Since $v$ is semi-standard, we have $v_{r,1}>v_{r-1,1}$.
Therefore, if we insert $v_{r-1,1}$, $v_{r,1}$ is bumped down
and come to the next to the row $22\cdots 2$.
Insertion of $v_{r-1,2}$ is quite similar, i.e.,
$v_{r,2}$ is bumped down and come to the next to
the row $22\cdots 2v_{r,1}$.
Continuing in this way, we can show that the resulting
tableau takes the following shape.
If $r\leq a$, we have
\begin{equation}
\left(u_k^{(a)}\leftarrow row(v)\right)=
\Yvcentermath1
\Yboxdim18pt
\newcommand{\vichiichi}{v_{1,1}}
\newcommand{\vichini}{v_{1,2}}
\newcommand{\vichis}{v_{1,s}}
\newcommand{\vniichi}{v_{2,1}}
\newcommand{\vnini}{v_{2,2}}
\newcommand{\vnis}{v_{2,s}}
\newcommand{\vrichi}{v_{r,1}}
\newcommand{\vrni}{v_{r,2}}
\newcommand{\vrs}{v_{r,s}}
\newcommand{\cdotsr}{\begin{picture}(2,2)
\multiput(-4,4)(4.7,0){3}{\circle*{1}}
\end{picture}}
\young(11\cdotsr 1\vichiichi\vichini\cdotsr\vichis,%
22\cdotsr 2\vniichi\vnini\cdotsr\vnis,%
\cdotsr\cdotsr\cdotsr\cdotsr\cdotsr\cdotsr\cdotsr\cdotsr,%
rrrr\vrichi\vrni\cdotsr\vrs,%
\cdotsr\cdotsr\cdotsr\cdotsr,%
aa\cdotsr a)\,\, ,
\end{equation}
and if $r>a$, we have
\begin{equation}
\left(u_k^{(a)}\leftarrow row(v)\right)=
\Yvcentermath1
\Yboxdim18pt
\newcommand{\vichiichi}{v_{1,1}}
\newcommand{\vichini}{v_{1,2}}
\newcommand{\vichis}{v_{1,s}}
\newcommand{\vniichi}{v_{2,1}}
\newcommand{\vnini}{v_{2,2}}
\newcommand{\vnis}{v_{2,s}}
\newcommand{\vaichi}{v_{a,1}}
\newcommand{\vani}{v_{a,2}}
\newcommand{\vas}{v_{a,s}}
\newcommand{\vrichi}{v_{r,1}}
\newcommand{\vrni}{v_{r,2}}
\newcommand{\vrs}{v_{r,s}}
\newcommand{\cdotsr}{\begin{picture}(2,2)
\multiput(-4,4)(4.7,0){3}{\circle*{1}}
\end{picture}}
\young(11\cdotsr 1\vichiichi\vichini\cdotsr\vichis,%
22\cdotsr 2\vniichi\vnini\cdotsr\vnis,%
\cdotsr\cdotsr\cdotsr\cdotsr\cdotsr\cdotsr\cdotsr\cdotsr,%
aaaa\vaichi\vani\cdotsr\vas,%
\cdotsr\cdotsr\cdotsr\cdotsr,%
\vrichi\vrni\cdotsr\vrs)\,\, .
\end{equation}
In both cases, the shape of
$(u_k^{(a)}\leftarrow row(v))$ coincides with
the concatenation of shapes of $u_k^{(a)}$ and $v$.
By Proposition \ref{prop:shimozono}
the corresponding energy is 0.
\hfill\rule{4pt}{10pt}

\begin{lemma}\label{lem:red_of_combR}
Consider the following four elements:
$$v=
\Yvcentermath1
\Yboxdim18pt
\newcommand{\cdotsr}{\begin{picture}(2,2)
\multiput(-4,4)(4.7,0){3}{\circle*{1}}
\end{picture}}
\newcommand{\bichi}{b_1}
\newcommand{\bni}{b_2}
\newcommand{\bs}{b_s}
\newcommand{\bpichi}{b'_1}
\newcommand{\bpni}{b'_2}
\newcommand{\bps}{b'_{s'}}
\newcommand{\tbpichi}{\tilde{b}'_1}
\newcommand{\tbpni}{\tilde{b}'_2}
\newcommand{\tbps}{\tilde{b}'_{s'}}
\newcommand{\tbichi}{\tilde{b}_1}
\newcommand{\tbni}{\tilde{b}_2}
\newcommand{\tbs}{\tilde{b}_s}
\young(11\cdotsr 1,22\cdotsr 2,%
\cdotsr\cdotsr\cdotsr\cdotsr,aa\cdotsr a,%
\bichi\bni\cdotsr \bs)
\, ,\,
v'=
\young(11\cdotsr 1,22\cdotsr 2,%
\cdotsr\cdotsr\cdotsr\cdotsr,aa\cdotsr a,%
\bpichi\bpni\cdotsr \bps)
\, ,\,
\tilde{v}'=
\young(11\cdotsr 1,22\cdotsr 2,%
\cdotsr\cdotsr\cdotsr\cdotsr,aa\cdotsr a,%
\tbpichi\tbpni\cdotsr \tbps)
\, ,\,
\tilde{v}=
\young(11\cdotsr 1,22\cdotsr 2,%
\cdotsr\cdotsr\cdotsr\cdotsr,aa\cdotsr a,%
\tbichi\tbni\cdotsr \tbs)
\, ,\,
$$
where $v,\tilde{v}\in B^{a+1,s}$ and $v',\tilde{v}'\in B^{a+1,s'}$.
The upper rows of $v,\tilde{v}$
$($resp. $v',\tilde{v}')$ coincide with
the highest elements $u_s^{(a)}\in B^{a,s}$
$($resp. $u_{s'}^{(a)}\in B^{a,s'})$.
The bottom rows satisfy
$b_i,b'_i,\tilde{b}_i,\tilde{b}'_i\in\{ a+1,a+2\}$.
The relationship between $b_i,b'_i,\tilde{b}_i,\tilde{b}'_i$
is given by
\begin{equation}
\Yvcentermath1
\Yboxdim18pt
\newcommand{\cdotsr}{\begin{picture}(2,2)
\multiput(-4,4)(4.7,0){3}{\circle*{1}}
\end{picture}}
\newcommand{\bichi}{b_1}
\newcommand{\bni}{b_2}
\newcommand{\bs}{b_s}
\newcommand{\bpichi}{b'_1}
\newcommand{\bpni}{b'_2}
\newcommand{\bps}{b'_{s'}}
\newcommand{\tbpichi}{\tilde{b}'_1}
\newcommand{\tbpni}{\tilde{b}'_2}
\newcommand{\tbps}{\tilde{b}'_{s'}}
\newcommand{\tbichi}{\tilde{b}_1}
\newcommand{\tbni}{\tilde{b}_2}
\newcommand{\tbs}{\tilde{b}_s}
\young(\bichi\bni\cdotsr \bs)
\otimes
\young(\bpichi\bpni\cdotsr \bps)
\simeq
\young(\tbpichi\tbpni\cdotsr \tbps)
\otimes
\young(\tbichi\tbni\cdotsr \tbs)
\, ,\,
\end{equation}
under the isomorphism $B^{1,s}\otimes B^{1,s'}\simeq
B^{1,s'}\otimes B^{1,s}$.

(1) If $k\neq a+1$, we have
\begin{equation}
u_l^{(k)}\otimes v\simeq v\otimes u_l^{(k)},\qquad
H\left(u_l^{(k)}\otimes v\right)=0.
\end{equation}

(2) We have
\begin{equation}
v\otimes v'\simeq\tilde{v}'\otimes \tilde{v},\qquad
H(v\otimes v')=
H\left(
\Yvcentermath1
\Yboxdim17pt
\newcommand{\cdotsr}{\begin{picture}(2,2)
\multiput(-4,4)(4.7,0){3}{\circle*{1}}
\end{picture}}
\newcommand{\bichi}{b_1}
\newcommand{\bni}{b_2}
\newcommand{\bs}{b_s}
\newcommand{\bpichi}{b'_1}
\newcommand{\bpni}{b'_2}
\newcommand{\bps}{b'_{s'}}
\young(\bichi\bni\cdotsr \bs)
\otimes
\young(\bpichi\bpni\cdotsr \bps)
\right).
\end{equation}
In other words, the combinatorial $R$ and
energy function for $v\otimes v'$ essentially coincide
with the $A_1^{(1)}$ type ones if we
replace $a+1$ and $a+2$ of $b_i, b'_i, \tilde{b}_i, \tilde{b}'_i$
with $1$ and $2$, respectively.
\end{lemma}
{\bf Proof.}
(1) We use the similar argument given in
Lemma \ref{lem:energy_conservation}.
In order to show the isomorphism of the combinatorial $R$,
one can directly compute $(v\leftarrow row(u_l^{(k)}))$
and $(u_l^{(k)}\leftarrow row(v))$ and
show that resulting two tableaux coincide.

(2) During the computation of $(v'\leftarrow row(v))$,
we arrive at the following intermediate diagram:
\begin{equation}
\Yvcentermath1
\Yboxdim17pt
\newcommand{\cdotsr}{\begin{picture}(2,2)
\multiput(-4,4)(4.7,0){3}{\circle*{1}}
\end{picture}}
\newcommand{\bichi}{b_1}
\newcommand{\bni}{b_2}
\newcommand{\bs}{b_s}
\newcommand{\bpichi}{b'_1}
\newcommand{\bpni}{b'_2}
\newcommand{\bps}{b'_{s'}}
\young(11\cdotsr 122\cdotsr 2,%
22\cdotsr 2\cdotsr\cdotsr\cdotsr\cdotsr,%
\cdotsr\cdotsr\cdotsr\cdotsr\cdotsr\cdotsr\cdotsr\cdotsr,%
\cdotsr\cdotsr\cdotsr\cdotsr aa\cdotsr a,%
aa\cdotsr a\bichi\bni\cdotsr\bs,%
\bpichi\bpni\cdotsr\bps) \, .
\end{equation}
As the final step of $(v'\leftarrow row(v))$,
we insert the row $11\cdots 1$ into this diagram.
If we insert 1, we see that $b_1$ is bumped,
and we have to calculate the insertion of $b_1$
into the bottom row $b_1'b_2'\cdots b'_{s'}$.
Continuing in this way, we see that the calculation
essentially coincides with that for
$(b_1'b_2'\cdots b'_{s'}\leftarrow
b_1b_2\cdots b_s)$,
hence we finish the proof.
\hfill\rule{4pt}{10pt}

\subsection{Reduction to the $A_1^{(1)}$ case}
We reduce the proof of Proposition \ref{prop:E=Q}
to the $A_1^{(1)}$ case by invoking analogous trick
described in Section 2.7 of \cite{KOSTY}.
These kinds of arguments are reformulation of
color separation scheme of the box-ball systems
\cite{Takagi} from point of view of the KSS bijection.
Given a path
$b=b_1\otimes b_2\otimes\cdots\otimes b_L
\in B^{\alpha_1,\beta_1}\otimes B^{\alpha_2,\beta_2}
\otimes\cdots \otimes B^{\alpha_L,\beta_L}$,
let the corresponding rigged configuration be
\begin{equation}
\mathrm{RC}=
\left(
(\nu_i^{(0)})_{i=1}^{L^{(0)}},\cdots,
(\nu_i^{(n-1)})_{i=1}^{L^{(n-1)}},
(\mu_i^{(1)},r_i^{(1)})_{i=1}^{N^{(1)}},\cdots,
(\mu_i^{(n)},r_i^{(n)})_{i=1}^{N^{(n)}}
\right).
\end{equation}
We consider the following modification of this RC:
\begin{eqnarray}
\mathrm{RC}_+&=&
\left(
(\nu_i^{(0)})_{i=1}^{L^{(0)}}\cup (1^{l^{(0)}+m^{(0)}}),\cdots,
(\nu_i^{(n-1)})_{i=1}^{L^{(n-1)}}\cup (1^{l^{(n-1)}+m^{(n-1)}}),
\right.\nonumber\\
&&\qquad\left.
(\mu_i^{(1)},r_i^{(1)}+l^{(0)})_{i=1}^{N^{(1)}},\cdots,
(\mu_i^{(n)},r_i^{(n)}+l^{(n-1)})_{i=1}^{N^{(n)}}
\right).
\end{eqnarray}
We take $l^{(a)}\gg |\mu^{(a)}|$ ($l^{(0)}$ is arbitrary), and
integers $m^{(a)}$ are chosen sufficiently large
according to the procedure described in the following.

We apply $\phi^{-1}$ on this RC${}_+$ in
the following two different ways.
First, we remove $(1^{m^{(n-1)}})$, $\cdots$,
$(1^{m^{(0)}})$ from the quantum space,
and then remove the rest of the quantum space $\nu^{(a)}$ by the
same way to obtain $b$ from RC.
This is possible since the change of the vacancy number
induced by the extra $(1^{l^{(a)}})$ is canceled by the
shift of riggings.
Finally, we are left with the rigged configuration
\begin{equation}
\left(
(1^{l^{(0)}}),\cdots,(1^{l^{(n-1)}}),
(\emptyset ,\emptyset),\cdots (\emptyset ,\emptyset)
\right)
\end{equation}
and the calculation of $\phi^{-1}$ on it becomes trivial.
Therefore, the resulting path $\tilde{b}$ takes the following form:
\begin{align}
&\tilde{b}=\\
&\newcommand{\tatedot}{\begin{picture}(2,2)
\multiput(1,0)(0,3.5){3}{\circle*{1}}
\end{picture}}
\Yvcentermath1
\left.\young(1)\right.^{\,\otimes l^{(0)}}
\otimes\left.\young(1,2)\right.^{\,\otimes l^{(1)}}
\otimes\cdots\otimes
\left.\young(1,2,\tatedot,n)\right.^{\,\otimes l^{(n-1)}}
\otimes b\otimes
\left.\young(1)\right.^{\,\otimes m^{(0)}}
\otimes\left.\young(1,2)\right.^{\,\otimes m^{(1)}}
\otimes\cdots\otimes
\left.\young(1,2,\tatedot,n)\right.^{\,\otimes m^{(n-1)}}.
\nonumber
\end{align}

Next, we remove whole
$(\nu_i^{(n-1)})_{i=1}^{L^{(n-1)}}\cup
(1^{l^{(n-1)}+m^{(n-1)}})$
before removing other part of the quantum space.
The integers $m^{(a)}$ are taken to be large enough
such that no rows of $\mu^{(i)}$ ($1\leq i\leq n-1$)
become singular until entire removal of
$(\nu_i^{(n-1)})\cup (1^{l^{(n-1)}+m^{(n-1)}})$.
The necessity of the addition $(1^{l^{(n-1)}})$
is clarified from the following observation.
Let us compare the actual situation on $\mathrm{RC}_+$
with the following
$A^{(1)}_1$ rigged configuration
\begin{equation}\label{eq:A11reduction}
\left(
(\nu_i^{(n-1)})_{i=1}^{L^{(n-1)}}\cup (1^{l^{(n-1)}+m^{(n-1)}}),
(\mu_i^{(n)},r_i^{(n)})_{i=1}^{N^{(n)}}
\right).
\end{equation}
Then the vacancy numbers for $\mathrm{RC}_+$
are increased by the presence
of $\mu^{(n-1)}$, and therefore their increment
from Eq.(\ref{eq:A11reduction}) are at
most $|\mu^{(n-1)}|$.
Recall that in $\mathrm{RC}_+$,
the riggings for $\mu^{(n)}$ are increased by $l^{(n-1)}$,
and also that we took $l^{(n-1)}\gg|\mu^{(n-1)}|$.
Thus, for each row $\mu^{(n)}_i$ of $\mathrm{RC}_+$,
the value of the corigging (see the end of the first
paragraph of Section \ref{subsec:KSS}) is
smaller than that for the corresponding row
of Eq.(\ref{eq:A11reduction}).
Note that the coriggings are always non-negative,
and if it is zero, then the row is singular
and can be removed by $\phi^{-1}$.
In other words, smaller coriggings means
much easier to be removed.
Since all boxes of
$(\mu_i^{(n)},r_i^{(n)})_{i=1}^{N^{(n)}}$
of Eq.(\ref{eq:A11reduction}) are removed with
$(\nu_i^{(n-1)})_{i=1}^{L^{(n-1)}}\cup (1^{l^{(n-1)}+m^{(n-1)}})$,
we deduce that whole of
$(\mu_i^{(n)},r_i^{(n)}+l^{(n-1)})_{i=1}^{N^{(n)}}$
of $\mathrm{RC}_+$
is removed  with
$(\nu_i^{(n-1)})_{i=1}^{L^{(n-1)}}\cup (1^{l^{(n-1)}+m^{(n-1)}})$.
As the result, we obtain the path
which we denote by $b^{(n)}$,
and the rest of the rigged configuration takes
the following form:
\begin{align}
&\left(
(\nu_i^{(0)})_{i=1}^{L^{(0)}}\cup (1^{l^{(0)}+m^{(0)}}),\cdots,
(\nu_i^{(n-2)})_{i=1}^{L^{(n-2)}}\cup (1^{l^{(n-2)}+m^{(n-2)}}),
\right.\\
&\qquad\left.(\mu_i^{(1)},r_i^{(1)}+l^{(0)})_{i=1}^{N^{(1)}},\cdots,
(\mu_i^{(n-1)},r_i^{(n-1)}+l^{(n-2)})_{i=1}^{N^{(n-1)}}
\right).
\end{align}
Since $\mu^{(a)}$ ($1\leq a\leq n-1$) do not
become singular during this procedure,
each tensor factor of $b^{(n)}$
have special property:
as a tableau, the upper $n-1$ rows are equal to
the highest element, and the bottom row
consists of letters $n$ and $n+1$.
In this sense, $b^{(n)}$
is obtained by the $A^{(1)}_1$ like KSS bijection.
Next, we remove whole
$(\nu_i^{(n-2)})_{i=1}^{L^{(n-2)}}\cup (1^{l^{(n-2)}+m^{(n-2)}})$
in a similar manner.
We denote the resulting path $b^{(n-1)}$.
We repeat the procedure until whole RC${}_+$ is removed.
Then, in general, tensor factors of $b^{(a+1)}$
have the following form:
\begin{equation}
\newcommand{\bichi}{b_1}
\newcommand{\bni}{b_2}
\newcommand{\bs}{b_s}
\Yvcentermath1
\Yboxdim17pt
\newcommand{\cdotsr}{\begin{picture}(2,2)
\multiput(-4,4)(4.7,0){3}{\circle*{1}}
\end{picture}}
\young(11\cdotsr 1,22\cdotsr 2,%
\cdotsr\cdotsr\cdotsr\cdotsr ,aa\cdotsr a,%
\bichi\bni\cdotsr\bs)\, ,
\end{equation}
where the bottom row consisting of
$b_i\in\{ a+1,a+2\}$.
Finally, we obtain the path
\begin{equation}
\tilde{\tilde{b}}=
b^{(1)}\otimes b^{(2)}\otimes\cdots\otimes
b^{(n)}.
\end{equation}
Since the difference between
these two ways of removal is order of removing
rows of the quantum space, these two paths are
isomorphic
\begin{equation}
\tilde{b}\simeq\tilde{\tilde{b}}
\end{equation}
due to Theorem \ref{th:KSS}.
\bigskip

\noindent
{\bf Proof of Proposition \ref{prop:E=Q}.}
{}From Lemma \ref{lem:energy_conservation},
we have $E_l^{(a)}(b)=E_l^{(a)}\bigl(\tilde{b}\bigl)$.
Since we have isomorphism
$\tilde{b}\simeq\tilde{\tilde{b}}$,
we have $E_l^{(a)}\bigl(\tilde{b}\bigl)=
E_l^{(a)}\bigl(\tilde{\tilde{b}}\bigl)$
thanks to the Yang--Baxter relation
(Proposition \ref{prop:YBeq}).
To see this, for the isomorphism
$x\otimes y\simeq\tilde{y}\otimes\tilde{x}$
with $H(x\otimes y)=e$, let us assign the diagram
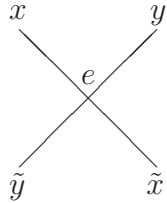
\unitlength 13pt
\begin{center}
\begin{picture}(5,5)(0,-0.5)
\put(0,0){\line(1,1){4}}
\put(4,0){\line(-1,1){4}}
\put(1.8,2.4){$e$}
\put(-0.3,4.3){$x$}
\put(3.8,4.3){$y$}
\put(-0.3,-0.8){$\tilde{y}$}
\put(3.7,-0.8){$\tilde{x}$}
\end{picture}
\end{center}
Then for tensor product $x\otimes y\otimes z$,
the Yang--Baxter relation takes
the following expression:
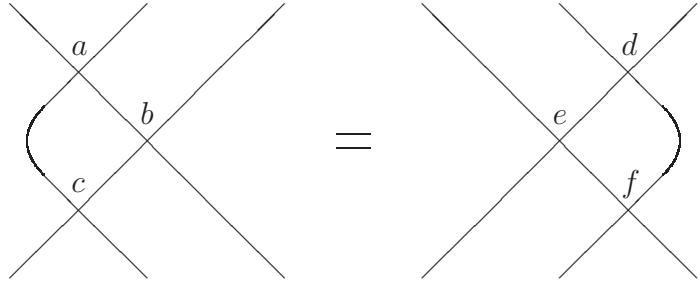
\begin{center}
\begin{picture}(20,8)
\put(0,0){\line(1,1){8}}
\put(8,0){\line(-1,1){8}}
\put(4,0){\line(-1,1){3}}
\put(4,8){\line(-1,-1){3}}
\qbezier(1,3)(0,4)(1,5)
\put(1.8,6.5){$a$}
\put(3.8,4.5){$b$}
\put(1.8,2.5){$c$}
\multiput(9.5,3.8)(0,0.4){2}{\line(1,0){1}}
\put(12,0){\line(1,1){8}}
\put(20,0){\line(-1,1){8}}
\put(16,0){\line(1,1){3}}
\put(16,8){\line(1,-1){3}}
\qbezier(19,3)(20,4)(19,5)
\put(17.8,6.5){$d$}
\put(15.8,4.5){$e$}
\put(17.8,2.5){$f$}
\end{picture}
\end{center}
In particular, we can show $a+b=e+f$
(see after Eq.(4.14)
of Section 4.3 of \cite{KSY} for more details).
If we put $x=u^{(a)}_l$, we obtain
$E^{(a)}_l(y\otimes z)=E^{(a)}_l(\tilde{z}\otimes\tilde{y})$
where $y\otimes z\simeq\tilde{z}\otimes\tilde{y}$.
More general case can be deduced inductively
from this case.

Now let us evaluate $E_l^{(a)}\bigl(\tilde{\tilde{b}}\bigl)$.
{}From Lemma \ref{lem:red_of_combR} (1),
we have $E_l^{(a)}\bigl(\tilde{\tilde{b}}\bigl)=
E_l^{(a)}\bigl( b^{(a)}\bigl)$.
Here, we use the following observation.
Since integer $m^{(a-1)}$ is large enough,
right part of $b^{(a)}$ consists of a large number of highest
elements.
Therefore, due to Lemma \ref{lem:red_of_combR} (2),
the carrier (see Eq.(\ref{def:T_l}) and
the comments following it) becomes
the highest element $u_l^{(a)}$,
hence we can ignore $b^{(a+1)}$, $\cdots$, $b^{(n)}$.
By construction, each $b^{(i)}$ is obtained from
$\mathrm{RC}_+$ by $\phi^{-1}$
on the quantum space
$(
(\nu_k^{(i-1)})_{k=1}^{L^{(i-1)}}\cup (1^{l^{(i-1)}+m^{(i-1)}}))$
and configuration
$((\mu_k^{(i)},r_k^{(i)}+l^{(i-1)})_{k=1}^{N^{(i)}})$.
If we compare this with the genuine
$A_1^{(1)}$ type KSS bijection on
\begin{equation}\label{eq:A11RC}
\left(
(\nu_k^{(i-1)})_{k=1}^{L^{(i-1)}}\cup (1^{l^{(i-1)}+m^{(i-1)}}),
(\mu_k^{(i)},r_k^{(i)}+l^{(i-1)})_{k=1}^{N^{(i)}}
\right),
\end{equation}
the difference is the change of
the vacancy numbers due to the presence of other $\mu^{(i)}$.
However these changes are constant during
construction of $b^{(i)}$.
Thus such differences can be absorbed into
the change of riggings.
In particular, in calculation of $Q_l^{(a)}$,
we can ignore such differences.
{}From this observation and Lemma \ref{lem:red_of_combR} (2),
we see that if we want to evaluate $E^{(a)}_l(b^{(a)})$,
it is enough to consider $A^{(1)}_1$ type rigged
configuration (\ref{eq:A11RC}) with $i=a$.

Therefore, we have to check Proposition \ref{prop:E=Q}
in $A_1^{(1)}$ case.
This can be done by a direct calculation
using time evolution of the box-ball systems.
Leaving precise arguments for Lemma 4.1 of \cite{Sak2},
we explain here main ideas of the proof.
Let $b\in B^{1,s_1}\otimes\cdots
\otimes B^{1,s_L}$ be an arbitraly tensor product of
$A^{(1)}_1$ crystals and consider the path
$b':=b\otimes \fbox{1}^{\,\otimes \Lambda}$,
where $\Lambda$ is sufficiently large integer.
Note that we have $E^{(1)}_l(b')=E^{(1)}_l(b)$
(see Eq.(\ref{eq:v*u=0})).
Consider time evolution $\bar{b}:=\left(T^{(1)}_s\right)^{t}(b')$
where $s$ and $t$ are sufficiently large integers.
{}From the inverse scattering method for the box-ball
systems (see Remark \ref{rem:ist}), we can show that
$\bar{b}$ has very simplified structure.
This enables us to compute $E^{(1)}_l(\bar{b})$ directly.
Again by use of the Yang--Baxter relation
(see Theorem 3.2 of \cite{FOY}),
we can show $E^{(1)}_l(\bar{b})=E^{(1)}_l(b')$.
Hence we complete the proof of Proposition \ref{prop:E=Q}.
\hfill\rule{4pt}{10pt}

\begin{remark}
Part of this section is based on joint work
with A. Kuniba (April 2007, unpublished).
\end{remark}

\section{Proof of Theorem \ref{th:main}}
To begin with, we have to settle some
combinatorial property
concerning the calculation of $\phi^{-1}$.
Given a path, we regard each semi-standard Young tableau
(i.e., tensor factor of the path)
as collection of columns.
We shall examine in details the construction of
these columns from the given rigged configuration.
Recall that $\phi^{-1}$ recursively constructs
tensor factors of the path from right to left,
and within each tensor factor, it recursively
constructs columns of tableau from left to right.
Therefore, without loss of generality,
we look at the leftmost column of the rightmost
tensor factor of the path.
We denote this column by $C$, and denote
the number of nodes of $C$ by $a+1$.
Then the construction of $C$ starts by
removing a box from $\nu^{(a)}$
(we name the box $x^{(a)}$).
In the box removing procedure of $\phi^{-1}$,
starting from the box $x^{(a)}$,
we remove other boxes (say $y\in\mu^{(a+1)}$,
$z\in\mu^{(a+2)}$, and so on) simultaneously.
Then we add a box to the first
column of $\nu^{(a-1)}$, which we call $x^{(a-1)}$.
The next step starts from $x^{(a-1)}$
and do the similar procedure.

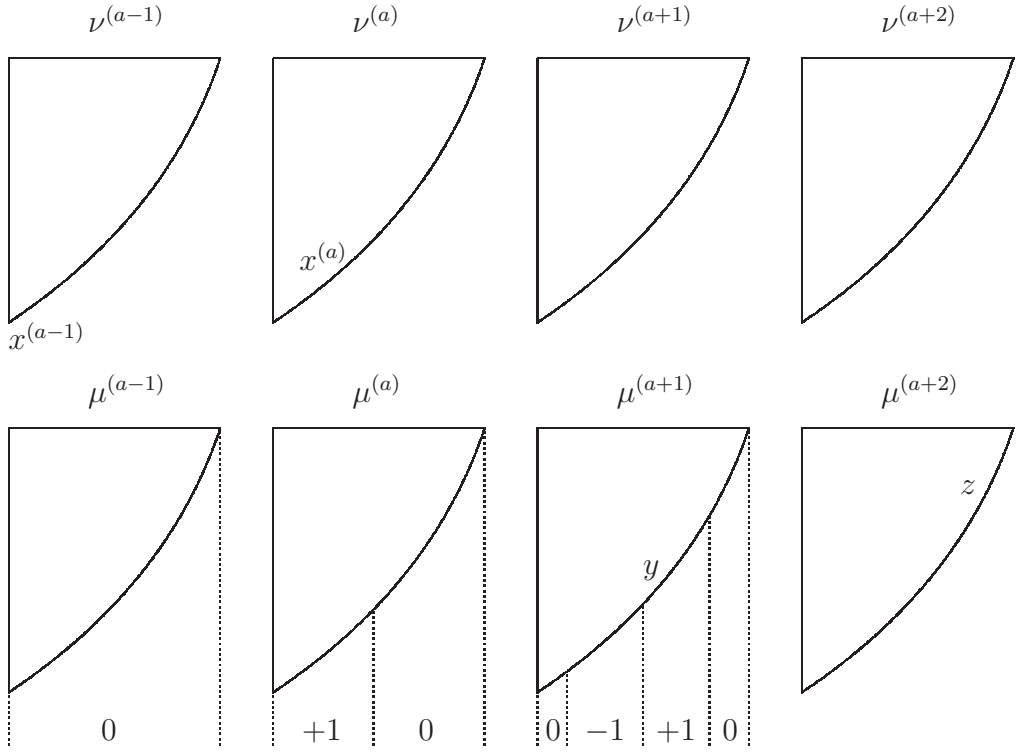
\begin{figure}
\begin{center}
\unitlength 10pt
\begin{picture}(40,28)(0,-2)
\multiput(0,0)(0,14){2}{
\multiput(0,0)(10,0){4}{
\put(0,0){\line(0,1){10}}
\put(0,10){\line(1,0){8}}
\qbezier(0,0)(6,4)(8,10)
}}
\put(3,11){$\mu^{(a-1)}$}
\put(13,11){$\mu^{(a)}$}
\put(23,11){$\mu^{(a+1)}$}
\put(33,11){$\mu^{(a+2)}$}
\put(3,25){$\nu^{(a-1)}$}
\put(13,25){$\nu^{(a)}$}
\put(23,25){$\nu^{(a+1)}$}
\put(33,25){$\nu^{(a+2)}$}
\put(0,13){$x^{(a-1)}$}
\put(11,16){$x^{(a)}$}
\put(24,4.5){$y$}
\put(36,7.5){$z$}
\multiput(0,-2)(0,0.2){10}{\line(0,1){0.07}}
\multiput(8,-2)(0,0.2){60}{\line(0,1){0.07}}
\put(3.5,-1.8){0}
\multiput(10,-2)(0,0.2){10}{\line(0,1){0.07}}
\multiput(13.8,-2)(0,0.2){26}{\line(0,1){0.07}}
\multiput(18,-2)(0,0.2){60}{\line(0,1){0.07}}
\put(11.1,-1.8){$+1$}
\put(15.5,-1.8){$0$}
\multiput(20,-2)(0,0.2){10}{\line(0,1){0.07}}
\multiput(21.1,-2)(0,0.2){14}{\line(0,1){0.07}}
\multiput(24.0,-2)(0,0.2){27}{\line(0,1){0.07}}
\multiput(26.5,-2)(0,0.2){44}{\line(0,1){0.07}}
\multiput(28,-2)(0,0.2){60}{\line(0,1){0.07}}
\put(20.3,-1.8){$0$}
\put(21.8,-1.8){$-1$}
\put(24.5,-1.8){$+1$}
\put(27,-1.8){$0$}
\end{picture}
\caption{Calculation of $\phi^{-1}$: remove $x^{(a)}$,
$y$, $z$, $\cdots$, and add $x^{(a-1)}$.}
\label{fig:phi}
\end{center}
\end{figure}

Let us summarize the situation in Figure \ref{fig:phi}.
Here, integers below the diagrams
show change of the vacancy numbers
caused by the removal of $x^{(a)}$,
$y$, $z$, $\cdots$, and the addition
of $x^{(a-1)}$.
We can determine these changes by
comparing number of boxes corresponding to
each term of definition Eq.(\ref{def:vacancy})
before and after the operation in the previous sentence.
Let us denote the column coordinate of a box $w$
in $\mu^{(i)}$ (or $\nu^{(i)}$) by $\mathrm{col}(w)$.
Then the vacancy number for
each row of $\mu^{(a)}$ whose rightmost box $w$ satisfies
$0\leq \mathrm{col}(w)<\mathrm{col}(y)$ increases
by 1, and those for the other region
do not change.
In the same way, the vacancy numbers for
the region
$\mathrm{col}(x^{(a)})\leq \mathrm{col}(w)
<\mathrm{col}(y)$ of $\mu^{(a+1)}$ decrease by 1,
those for the region
$\mathrm{col}(y)\leq \mathrm{col}(w)
<\mathrm{col}(z)$ are increased by 1,
and those for the other region do not change.

In the following lemma, we will
use these settings such as
the column $C$ is the
leftmost column of the rightmost
tensor factor of the path.
In other words, the column $C$ is the first
column constructed by the procedure $\phi^{-1}$.
\begin{lemma}\label{lem:evaluation_of_rigging}
Let us exclusively consider
the procedure $\phi^{-1}$ during
construction of the column $C$.
Assume that boxes of some row, say $\mu^i$
of $\mu^{(i)}$ $(1\leq i\leq n)$
are removed during construction of the column $C$.
Then the vacancy number for the row $\mu^i$
evaluated when its rightmost box (i.e., the firstly
removed box of the row) is removed
is equal to the one evaluated
before the construction of $C$.
\end{lemma}
{\bf Proof.}
Construction of the column $C$ under $\phi^{-1}$
begins with removal of the box $x^{(a)}$.
Assume that we have removed boxes $x^{(a)}$, $y$, $z$, $\cdots$
and added $x^{(a-1)}$.
See Figure \ref{fig:phi},
especially change of the vacancy numbers summarized there.
Let us examine which boxes can be removed from $\mu^{(a)}$,
$\mu^{(a+1)}$, $\cdots$ if we remove $x^{(a-1)}$.
First we consider $\mu^{(a)}$.
In $\mu^{(a)}$, there is no singular row within
the region whose vacancy numbers are increased
by removal of $x^{(a)}$, $y$, $z$, $\cdots$
and addition of $x^{(a-1)}$.
Therefore rows $w$ of $\mu^{(a)}$ that satisfy
$\mathrm{col}(y)\leq \mathrm{col}(w)$
are the candidate for removal.
We see that the vacancy numbers for
the rows $w$ of $\mu^{(a)}$ satisfying
$\mathrm{col}(y)\leq \mathrm{col}(w)$
do not change from those evaluated
before removing $x^{(a)}$.

Assume that we can remove some row of $\mu^{(a)}$,
then the next removal, if possible, will be a row of $\mu^{(a+1)}$.
Recall that by definition of $\phi^{-1}$, column coordinates of
simultaneously removed boxes are non decreasing
integer sequence (from left to right).
Thus we see that the next candidates for removal are within the region
of $\mu^{(a+1)}$ whose column coordinates are greater than
or equal to $\mathrm{col}(y)$.
Here we used the fact that we have removed a box of $\mu^{(a)}$
whose column coordinate is greater than
or equal to $\mathrm{col}(y)$ (see the end of the previous
paragraph).
Note that the row to which the box $y$ was originally belonged
is now shorter than $\mathrm{col}(y)$ since its length is
shortened by 1 due to the removal of the box $y$.
Note also that the vacancy numbers for rows $w$ of $\mu^{(a+1)}$
that satisfy
$\mathrm{col}(y)\leq \mathrm{col}(w)
<\mathrm{col}(z)$ are increased by 1
by removal of $x^{(a)}$, $y$, $z$, $\cdots$
and addition of $x^{(a-1)}$.
Therefore we see that candidates for removal
exist in the region of $\mu^{(a+1)}$ whose column coordinates
are greater than or equal to $\mathrm{col}(z)$.
In conclusion, the vacancy number for the row of $\mu^{(a+1)}$
that can be removed simultaneously with $x^{(a-1)}$ remains
the same as that evaluated before removing $x^{(a)}$.

We can continue the similar arguments for $\mu^{(a+2)}$,
$\mu^{(a+3)}$, $\cdots$.
To state the result for general case, let us consider $\mu^{(A)}$
for some $A>a+1$.
We divide the argument into the following two cases.
\begin{enumerate}
\item
Assume that some row of $\mu^{(A+1)}$ was removed with $x^{(a)}$.
Let us call this removed box of $\mu^{(A+1)}$ by $X$.
Then the situation essentially coincides with the
case of $\mu^{(a+1)}$ in place of $\mu^{(A)}$.
Namely, if we remove $x^{(a-1)}$, we can simultaneously remove, if possible,
a row of $\mu^{(A)}$ whose column coordinate is greater
than or equal to $\mathrm{col}(X)$.
In such case, the vacancy number for the row of $\mu^{(A)}$
that is removed with $x^{(a-1)}$ is the same as that
evaluated before removing $x^{(a)}$.

\item
On the other hand, if boxes of $\mu^{(A+1)}$
were not removed with $x^{(a)}$, then we see that
the boxes of $\mu^{(A)}$ will not be removed with
$x^{(a-1)}$.
To see this, let $A$ be the minimal integer
satisfying $A\geq a$ such that
the boxes of $\mu^{(A+1)}$ were not removed with $x^{(a)}$.
If $A=a$, then we remove only one box $x^{(a)}$ from $\nu^{(a)}$
and add a box $x^{(a-1)}$ at the first column of $\nu^{(a-1)}$.
Thus the vacancy numbers for all rows of $\mu^{(a)}$
are increased by 1.
Therefore we can not remove a box from $\mu^{(a)}$
simultaneously with removal of $x^{(a-1)}$.
In the following, we assume that $A>a$.
Let us denote the removed box of $\mu^{(A)}$ by $X$.
\begin{enumerate}
\item
After removing $x^{(a)}$, $y$, $z$, $\cdots$, $X$,
the vacancy numbers for
the rows of $\mu^{(A)}$ longer than or equal to
$\mathrm{col}(X)$ are increased by 1.
Therefore we can not remove rows of
$\mu^{(A)}$ longer than or equal to
$\mathrm{col}(X)$.
\item
On the other hand, we can apply the above case 1 to
$\mu^{(A-1)}$ to see that, if we remove a row of
$\mu^{(A-1)}$ with $x^{(a-1)}$, then the
removed row is longer than or equal to
$\mathrm{col}(X)$.
Therefore if we remove a box from $\mu^{(A)}$
with $x^{(a-1)}$, we have to remove
a row that is longer than or equal to
$\mathrm{col}(X)$.
Recall that by definition of $\phi^{-1}$,
column coordinates of
simultaneously removed boxes has to be
non-decreasing sequence from left to right.
\end{enumerate}
Combining these two mutually contradicting
conditions, we conclude
that we can not remove a box from $\mu^{(A)}$
with $x^{(a-1)}$.
This is nothing but the column strict semi-standard property
of the column $C$.
\end{enumerate}
To summarize, the vacancy numbers for all rows
that are removed simultaneously with $x^{(a-1)}$
are the same with those evaluated before
we remove $x^{(a)}$.

More general case corresponding to all $x^{(i)}$
can be shown inductively.
Namely, we replace $x^{(a)}$ and $x^{(a-1)}$
with $x^{(i)}$ and $x^{(i-1)}$, respectively,
to show the statement for $x^{(i-1)}$ case
from the assumption for $x^{(i)}$.
\hfill\rule{4pt}{10pt}

\bigskip

In order to connect the original procedure of $\phi^{-1}$
with the local energy
distribution, we introduce a table which
record the process of $\phi^{-1}$.
More precisely, we prepare array of $n$ tables
whose $a$-th table records the history of $\mu^{(a)}$
as follows.
Recall that in the calculation of $\phi^{-1}$,
we remove a box from the quantum space, say $\nu^{(i)}$ and,
at the same time, remove several boxes from $\mu^{(j)}$'s.
Then we add a box to $\nu^{(i-1)}$ and do the similar
procedure.
We do these procedures for several times (in this case,
$i+1$ times), then we obtain one column of a semi-standard
Young tableau of path.
Every time we add a column (say $C$) to tableau of the path,
we add one column to the $a$-th table (from
right to left according to $\phi^{-1}$).
The column in this $a$-th table contains the
following information: if we remove total of $m$ boxes
at $k$-th column of $\mu^{(a)}$ during the construction
of $C$, we put $m$ at the $k$-th row of the column
in the $a$-th table.
Furthermore, we draw curves  on the table which
join entries of the table row by row
if they belong to the same row of $\mu^{(a)}$.
Thus if an entry is $m$, then there are $m$
curves run through the position of $m$.

As a result, we obtain a table whose entries
are joined by curves.
Then, in general, these curves have
crossing with each other:
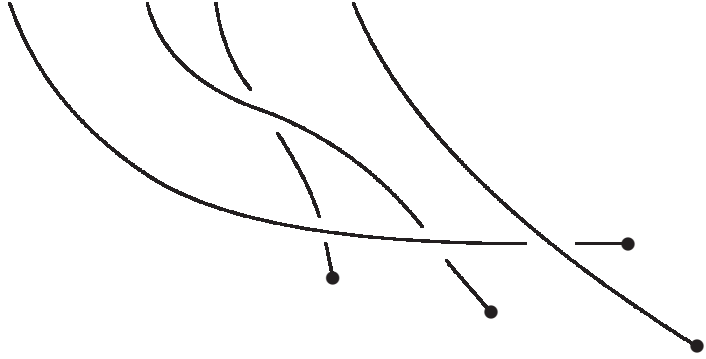
\begin{center}
\unitlength 13pt
\begin{picture}(20,10)
\thicklines
\qbezier(0,10)(1,7)(4,5)
\qbezier(4,5)(7,3)(15,3)
\qbezier(16.5,3)(18,3)(18,3)
\put(18,3){\circle*{0.4}}
\qbezier(4,10)(4.5,8)(7,7)
\qbezier(7,7)(10,6)(12,3.5)
\qbezier(12.7,2.5)(14,1)(14,1)
\put(14,1){\circle*{0.4}}
\qbezier(6,10)(6.2,8.5)(7,7.5)
\qbezier(7.8,6.2)(8.7,4.8)(9,3.8)
\qbezier(9.2,3)(9.2,3)(9.4,2)
\put(9.4,2){\circle*{0.4}}
\qbezier(10,10)(12,5)(20,0)
\put(20,0){\circle*{0.4}}
\end{picture}
\end{center}
Here parities of crossings (i.e., which line crosses over
the other line) are
not important.
Note that by construction of the curves joining
entries of table, these curves always go rightwards
and/or downwards, and never go leftwards
and/or upwards when we start from the top end.
In order to detect each lower end point $\bullet$,
we start from the rightmost strictly positive
integer in the top row.
We always go to weakly right of lower rows
(do not go to left),
and choose positive
integers row by row.
If we encounter the crossing,
we resolve it as follows:
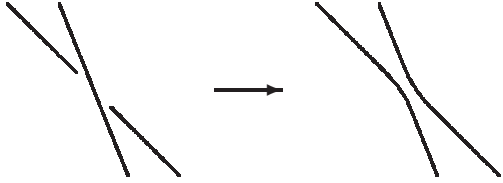
\begin{center}
\unitlength 13pt
\begin{picture}(14,5)
\thicklines
\qbezier(0,5)(0,5)(2,3)
\qbezier(3,2)(3,2)(5,0)
\qbezier(2.3,3)(2.3,3)(2.7,2)
\qbezier(1.5,5)(1.5,5)(2.3,3)
\qbezier(2.7,2)(3.5,0)(3.5,0)
\put(6,2.5){\vector(1,0){2}}
\put(9,0){
\qbezier(0,5)(2,3)(2,3)
\qbezier(2,3)(2.5,2.5)(2.7,2)%
\qbezier(3.3,2)(3.3,2)(5.3,0)
\qbezier(3.3,2)(2.9,2.4)(2.6,3)%
\qbezier(1.8,5)(1.8,5)(2.6,3)
\qbezier(2.7,2)(3.5,0)(3.5,0)
}
\end{picture}
\end{center}
In other words, if we have more than two
non-zero positive integers in the lower row
which are located weakly right to the chosen integer
of the upper row,
then we always choose the rightmost one, and proceed.
In this way, we arrive at one of the end points.
We subtract 1 from the chosen positive integers,
and choose the next group.
Continuing in this way, we can always
detect all the end points.

\bigskip

\noindent
{\bf Proof of Theorem \ref{th:main}.}
First, let us show the equivalence between the local energy
distribution and the above table recording the process
of $\phi^{-1}$.
{}From Proposition \ref{prop:E=Q},
we give interpretation of differences of
$E_{l,j,k}^{(a)}$
of Eq.(\ref{def:E_ljk}) as follows.
As before, consider the path
$b=b_1\otimes b_2\otimes\cdots\otimes b_L\in
B^{\alpha_1,\beta_1}\otimes B^{\alpha_2,\beta_2}\otimes
\cdots\otimes B^{\alpha_L,\beta_L}$,
and denote the columns of tableau $b_j$
as $b_j=c_{\beta_j}\cdots c_2c_1$,
and set $b_{j,k}=c_k\cdots c_2c_1$.
Consider the path
$b_{[k]}=b_1\otimes b_2\otimes\cdots\otimes b_{j-1}
\otimes b_{j,k}$
and apply Proposition \ref{prop:E=Q}.
Then we see that
\begin{equation}
\sum_{i=1}^{j-1}E_{l,i,\beta_i}^{(a)}+E_{l,j,k}^{(a)}
=(Q_l^{(a)}\mbox{ for the path }b_{[k]}).
\end{equation}
Therefore we have
\begin{equation}
E_{l,j,k}^{(a)}-E_{l,j,k-1}^{(a)}
=(Q_l^{(a)}\mbox{ for the path }b_{[k]})-
(Q_l^{(a)}\mbox{ for the path }b_{[k-1]}).
\end{equation}
By comparing with the procedure of $\phi^{-1}$,
the meaning for the right hand side is clear.
Namely, it represents the number of boxes
removed from the first $l$ columns of $\mu^{(a)}$
during construction of the column $c_k$ under $\phi^{-1}$.
Note that $Q_l^{(a)}$ represents, by definition,
the number of boxes
contained in the first $l$ columns of $\mu^{(a)}$.
Therefore taking difference once more,
\begin{equation}
\left(E_{l,j,k}^{(a)}-E_{l,j,k-1}^{(a)}\right)-
\left(E_{l-1,j,k}^{(a)}-E_{l-1,j,k-1}^{(a)}\right),
\end{equation}
it represents the number of boxes removed
from $l$-th column of $\mu^{(a)}$
corresponding to the column $c_k$
of a path under $\phi^{-1}$.
Therefore, the table introduced just before
the proof is equivalent to the
local energy distribution.
{}From this fact, it also follows that entries
of the local energy distribution are non-negative.
Furthermore, the end points detected in the
arguments before the proof correspond
to the rightmost box of each row of partitions $\mu^{(a)}$.
This follows from the fact that process of
$\phi^{-1}$ proceeds from right to left
in the table.
Hence we see that Steps 2 and 3 of the procedure given
in Theorem \ref{th:main} certainly extract the
data $\mu^{(a)}$ $(1\leq a\leq n)$.

Finally, it remains to show that
Eq.(\ref{eq:rigging1}) to Eq.(\ref{eq:rigging3})
in Step 4 really give the riggings.
Consider the arbitrary end point appearing
in the local energy distribution,
say at $\mu_i^{(a)}$-th row and $(j,k)$-th column.
As before, consider the path
$b_1\otimes b_2\otimes\cdots\otimes b_{j,k}$.
Then the end point we are considering
means that the corresponding row $\mu_i^{(a)}$
is removed for the first time
when we are constructing the column $c_k$
of the factor $b_j$ under $\phi^{-1}$.
By definition of $\phi^{-1}$,
the rigging is equal to the corresponding
vacancy number when the row $\mu_i^{(a)}$
is removed, and by
Lemma \ref{lem:evaluation_of_rigging},
it is enough to evaluate the vacancy number
at the time when we begin to construct $c_k$.
This means that it is enough to consider the path
$b_1\otimes b_2\otimes\cdots\otimes b_{j,k}$,
and evaluate the vacancy number for
the corresponding rigged configuration.
Then $\mathcal{C}$ of Eq.(\ref{eq:rigging2})
counts the contribution from the quantum space
(the first term of the right hand side of Eq.(\ref{def:vacancy})).
$\delta_{\alpha_i,a}$ in the expression Eq.(\ref{eq:rigging2})
reflects the
fact that only the contribution from $\nu^{(a-1)}$
appear in the vacancy number,
and the second term reflects the fact that we are
considering the truncated path
$b_1\otimes b_2\otimes\cdots\otimes b_{j,k}$.
The $\mathcal{E}$ of Eq.(\ref{eq:rigging3})
counts the contribution from $\mu^{(a-1)}$,
$\mu^{(a)}$ and $\mu^{(a+1)}$ of the
truncated path
(the last three terms of the
right hand side of Eq.(\ref{def:vacancy})).

We see that Steps 1 to 5 appearing
in the procedure described by
our theorem permit to recover
the rigged configuration we
are considering. This proves the theorem.
\hfill\rule{4pt}{10pt}

\appendix
\section{Kirillov--Schilling--Shimozono bijection}
\subsection{Rigged configurations}
In this appendix, we collect necessary facts about the
Kirillov--Schilling--Shimozono (KSS) bijection \cite{KSS}.
The reader should consult review \cite{Sch1}
for more details.
The KSS bijection gives one to one correspondences between
elements of tensor products $B^{r_1,s_1}\otimes
B^{r_2,s_2}\otimes\cdots\otimes B^{r_L,s_L}$ (which
we call {\it paths})
and combinatorial objects called the
rigged configurations.
The original theory deals with the highest weight elements,
and we first consider such situation.
Let us define the $A^{(1)}_n$ rigged configurations.
Consider the following set of data:
\begin{equation}
\mathrm{RC}=
\left(
(\nu_i^{(0)})_{i=1}^{L^{(0)}},\cdots,
(\nu_i^{(n-1)})_{i=1}^{L^{(n-1)}},
(\mu_i^{(1)},r_i^{(1)})_{i=1}^{N^{(1)}},\cdots,
(\mu_i^{(n)},r_i^{(n)})_{i=1}^{N^{(n)}}
\right).
\end{equation}
Here, $\nu_i^{(a)}$ and $\mu_i^{(a)}$ are
positive integer sequences.
Now we describe the conditions imposed on this RC.
Denote the number of boxes contained in
the first $l$ columns of the Young
diagrammatic expression of $\mu^{(a)}$
by $Q_l^{(a)}$, i.e.,
\begin{equation}\label{def:Q}
Q_l^{(a)}=\sum_{i=1}^{N^{(a)}}\min (l,\mu^{(a)}_i).
\end{equation}
Then, corresponding to the length $l$ row of
$\mu^{(a)}$, we define the {\it vacancy number}
by the following formula:
\begin{equation}\label{def:vacancy}
p_l^{(a)}=\sum_{i=1}^{L^{(a-1)}}\min (l,\nu_i^{(a-1)})
+Q_l^{(a-1)}-2Q_l^{(a)}+Q_l^{(a+1)},
\end{equation}
where the first term in the right hand side gives
the number of boxes contained in the first $l$
columns of $\nu^{(a-1)}$,
and we set $Q_l^{(0)}=Q_l^{(n+1)}=0$.
The integer $r_i^{(a)}$ associated with $\mu_i^{(a)}$
is called the {\it rigging associated with} $\mu_i^{(a)}$
if it satisfies the condition
\begin{equation}
0\leq r_i^{(a)}\leq p^{(a)}_{\mu_i^{(a)}}.
\end{equation}

Then the above RC is called the {\it rigged configuration}
if all the vacancy numbers are non-negative:
\begin{equation}
0\leq p^{(a)}_{\mu_i^{(a)}},\qquad
(1\leq a\leq n,\, 1\leq i\leq N^{(a)}),
\end{equation}
and if all integers $r_i^{(a)}$
are the riggings associated with $\mu_i^{(a)}$.
We usually call $\nu^{(a)}$ of RC the {\it quantum space},
and $\mu^{(a)}$ of RC the {\it configuration}.

\subsection{Combinatorial procedures for KSS bijection}\label{subsec:KSS}
By the KSS bijection, each path $b$ is mapped to
the rigged configuration RC:
\begin{equation}
\phi :b\longmapsto \mathrm{RC}.
\end{equation}
Below we describe how to calculate $\phi^{-1}$
rather than $\phi$.
The calculation of $\phi^{-1}$ consisting of
array of recursive procedures,
which are described by box removing processes.
In order to describe which box of $\mu^{(a)}$
to be removed, we use the following notion:
if the row $\mu_i^{(a)}$ satisfies the condition
\begin{equation}
p_{\mu_i^{(a)}}^{(a)}=r_i^{(a)},
\end{equation}
i.e., the rigging and the corresponding vacancy number
coincide, then the row $\mu_i^{(a)}$ is
called {\it singular}.
The value $p_{\mu_i^{(a)}}^{(a)}-r_i^{(a)}$
is called {\it corigging}.

The calculation of $\phi^{-1}$ starts by choosing
which row of $\nu^{(a)}$ $(0\leq a\leq n-1)$
we are going to remove.
Hence the map $\phi^{-1}$ is a function of such choice
(see Theorem \ref{th:KSS} below).
We choose row, say $\nu_i^{(a)}$, to be removed.
{}From the row $\nu_i^{(a)}$,
we obtain an element of the crystal
$B^{a+1,\nu_i^{(a)}}$ by the following procedure.

\begin{enumerate}
\item
Let us denote the rightmost box of the row
$\nu_i^{(a)}$ by $x^{(a)}$ (see Figure \ref{fig:phi}
in the main text).
Starting from $x^{(a)}$, we choose rows of
$\mu^{(a+1)}$, $\mu^{(a+2)}$, $\cdots$,
by the following recursive rule.
Assume that we have chosen a row of $\mu^{(j)}$.
Then we choose
a row of $\mu^{(j+1)}$ which is the shortest singular
row among the rows of $\mu^{(j+1)}$ not shorter
than the chosen row of $\mu^{(j)}$,
and proceed to $\mu^{(j+2)}$.
If there are more than one such row, choose one arbitrarily.
If there is not such a row, then stop.
If we can choose a row from $\mu^{(j_a-1)}$, and
cannot choose in $\mu^{(j_a)}$,
we obtain integer $j_a$ as the output.

\item
We remove all the rightmost boxes of chosen rows
simultaneously, and add a length 1 row
to $\nu^{(a-1)}$ which we call $x^{(a-1)}$.
After removal and addition,
we calculate the new vacancy numbers.
New riggings are determined as follows.
If the corresponding row is not removed,
choose the same rigging as before.
On the other hand, if the corresponding row
is removed, then choose the new rigging
equal to the new vacancy number for the
corresponding row.

\item
Repeat the above Step 1 and Step 2
starting from $x^{(a-1)}$, and obtain
$j_{a-1}$ as an output.
We do this until $x^{(0)}$ is removed.
Then we put $j_a$, $\cdots$, $j_{1}$,
$j_0$ to the leftmost empty column of
$(a+1)\times\nu_i^{(a)}$ Young diagram as follows:
\unitlength 1pt
\begin{equation}
\Yvcentermath1
\newcommand{\jeia}{j_a}
\newcommand{\jeiichi}{j_{1}}
\newcommand{\tatedot}{\begin{picture}(2,2)
\multiput(1,0)(0,3.5){3}{\circle*{1}}
\end{picture}}
\newcommand{\jeizero}{j_0}
\young(\jeizero\hfill\hfill\hfill\hfill\hfill,%
\jeiichi\hfill\hfill\hfill\hfill\hfill,%
\tatedot\hfill\hfill\hfill\hfill\hfill,%
\jeia\hfill\hfill\hfill\hfill\hfill)
\end{equation}

\item
Repeat the above Step 1 through Step 3
for the rest of $\nu_i^{(a)}$ box by box
(from right to left)
until entire row is removed.
Then we obtain an element of the crystal
$B^{a+1,\nu_i^{(a)}}$.
\end{enumerate}
We repeat the above Step 1 through Step 4
by choosing rows of the quantum space arbitrarily
until all boxes of the quantum space are removed.
We take tensor products of these tableaux from right
to left.

It is known that the above procedure is
well-defined, and the inverse map $\phi$
admits similar description \cite{KSS}.
In fact, the procedure for $\phi$
is obtained by just reversing the above procedure.

\begin{example}\label{ex:KSS}
Consider the following rigged configuration:
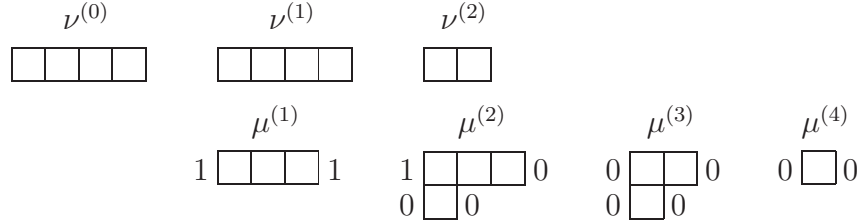
\begin{center}
\unitlength 13pt
\begin{picture}(26,6)
\Yboxdim13pt
\put(1,4){\young(\hfill\hfill\hfill\hfill)}
\put(2.5,5.5){$\nu^{(0)}$}
\put(7,1){\young(\hfill\hfill\hfill)}
\put(6.3,1.1){1}
\put(10.2,1.1){1}
\put(8,2.5){$\mu^{(1)}$}
\put(7,4){\young(\hfill\hfill\hfill\hfill)}
\put(8.5,5.5){$\nu^{(1)}$}
\put(13,0){\young(\hfill\hfill\hfill,\hfill)}
\put(12.3,0.1){0}
\put(12.3,1.1){1}
\put(14.2,0.1){0}
\put(16.2,1.1){0}
\put(14,2.5){$\mu^{(2)}$}
\put(13,4){\young(\hfill\hfill)}
\put(13.5,5.5){$\nu^{(2)}$}
\put(19,0){\young(\hfill\hfill,\hfill)}
\put(18.3,0.1){0}
\put(18.3,1.1){0}
\put(20.2,0.1){0}
\put(21.2,1.1){0}
\put(19.5,2.5){$\mu^{(3)}$}
\put(24,1){\young(\hfill)}
\put(23.3,1.1){0}
\put(25.2,1.1){0}
\put(24,2.5){$\mu^{(4)}$}
\end{picture}
\end{center}
Here we put the riggings (resp. vacancy numbers)
on the right (resp. left) of the corresponding rows.
As the first step, let us remove $\nu^{(1)}$.
Then calculation goes as follows
(we put $\times$ to the box to be removed):
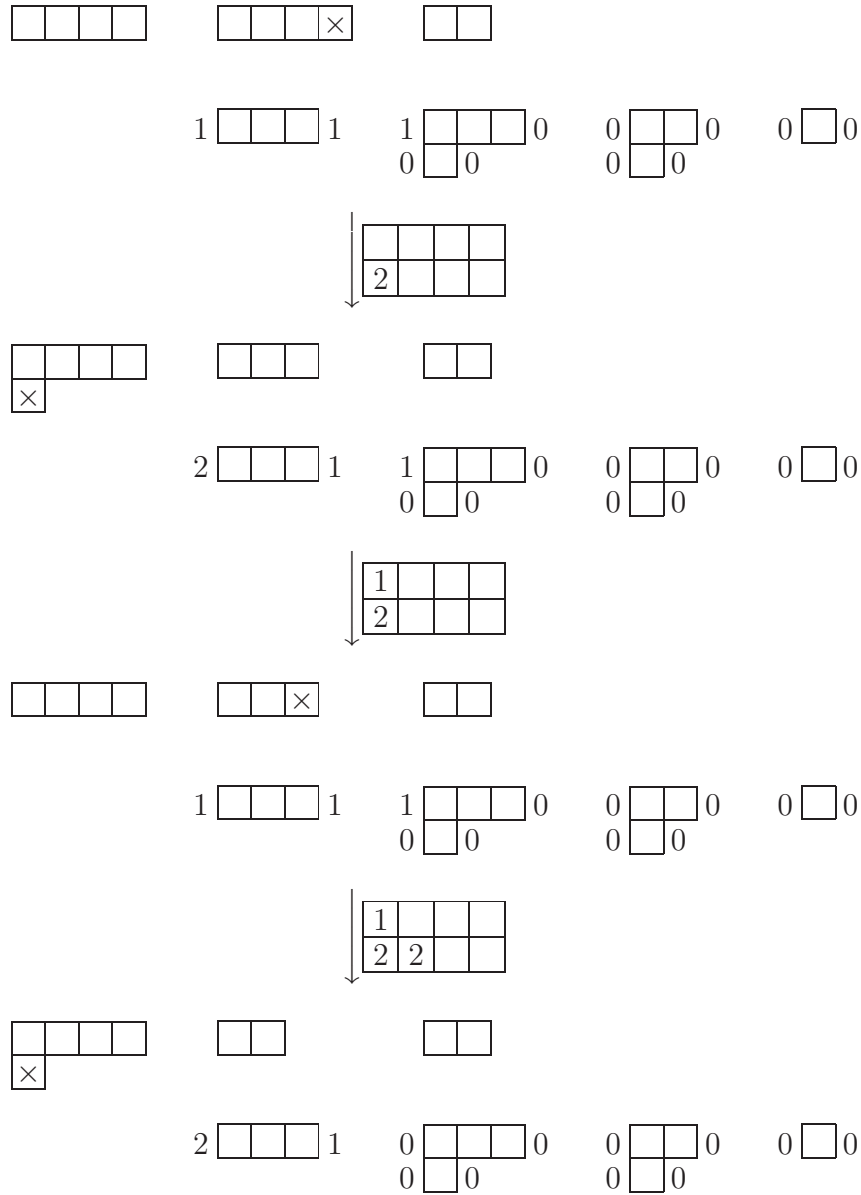
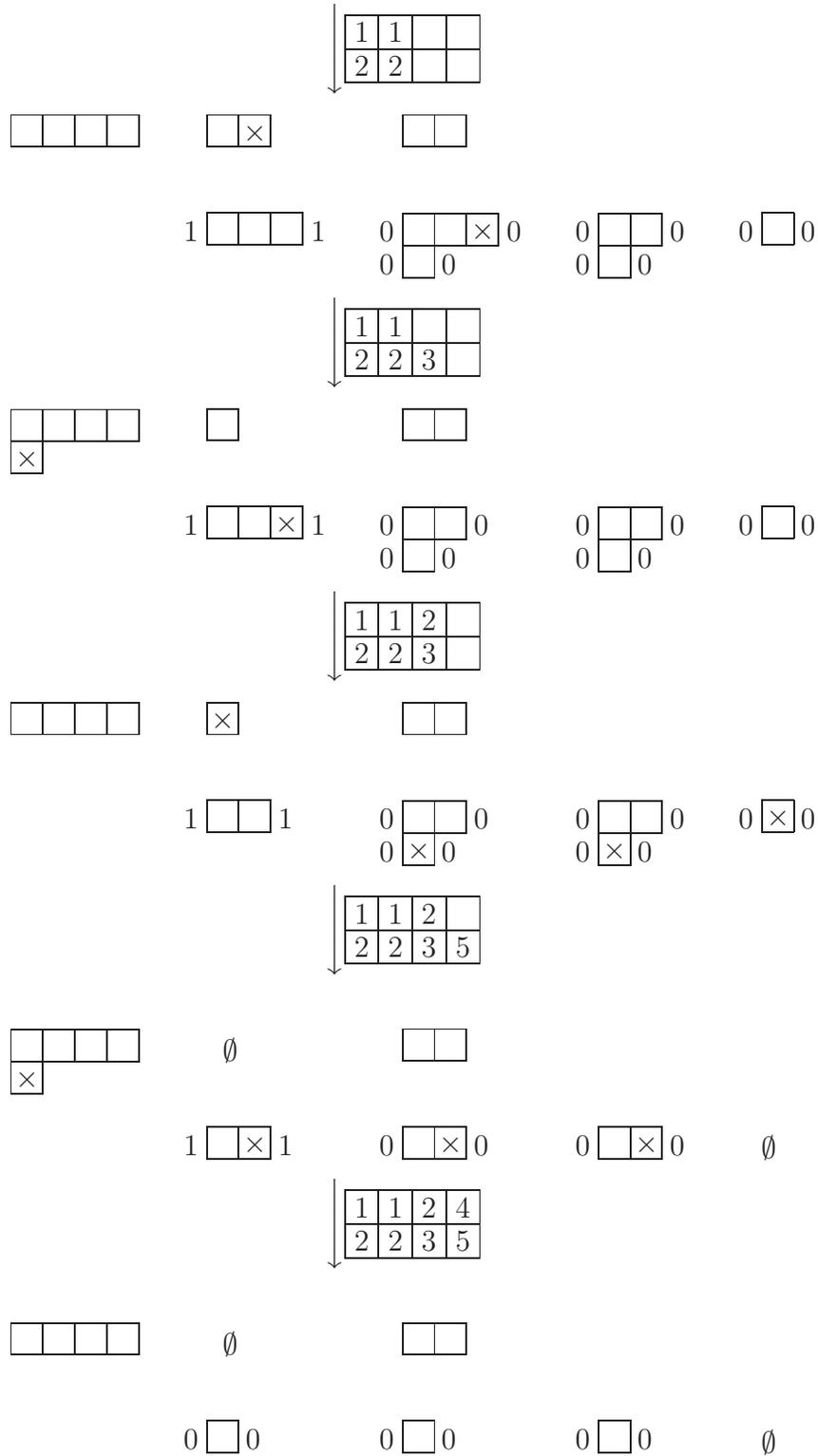
\begin{center}
\unitlength 13pt
\begin{picture}(26,5)
\Yboxdim13pt
\put(1,4){\young(\hfill\hfill\hfill\hfill)}
\put(7,1){\young(\hfill\hfill\hfill)}
\put(6.3,1.1){1}
\put(10.2,1.1){1}
\put(7,4){\young(\hfill\hfill\hfill\times)}
\put(13,0){\young(\hfill\hfill\hfill,\hfill)}
\put(12.3,0.1){0}
\put(12.3,1.1){1}
\put(14.2,0.1){0}
\put(16.2,1.1){0}
\put(13,4){\young(\hfill\hfill)}
\put(19,0){\young(\hfill\hfill,\hfill)}
\put(18.3,0.1){0}
\put(18.3,1.1){0}
\put(20.2,0.1){0}
\put(21.2,1.1){0}
\put(24,1){\young(\hfill)}
\put(23.3,1.1){0}
\put(25.2,1.1){0}
\end{picture}
\end{center}
\begin{center}
$\Bigg\downarrow
\Yvcentermath1
\young(\hfill\hfill\hfill\hfill,2\hfill\hfill\hfill)$
\end{center}
\begin{center}
\unitlength 13pt
\begin{picture}(26,5)
\Yboxdim13pt
\put(1,3){\young(\hfill\hfill\hfill\hfill,\times)}
\put(7,1){\young(\hfill\hfill\hfill)}
\put(6.3,1.1){2}
\put(10.2,1.1){1}
\put(7,4){\young(\hfill\hfill\hfill)}
\put(13,0){\young(\hfill\hfill\hfill,\hfill)}
\put(12.3,0.1){0}
\put(12.3,1.1){1}
\put(14.2,0.1){0}
\put(16.2,1.1){0}
\put(13,4){\young(\hfill\hfill)}
\put(19,0){\young(\hfill\hfill,\hfill)}
\put(18.3,0.1){0}
\put(18.3,1.1){0}
\put(20.2,0.1){0}
\put(21.2,1.1){0}
\put(24,1){\young(\hfill)}
\put(23.3,1.1){0}
\put(25.2,1.1){0}
\end{picture}
\end{center}
\begin{center}
$\Bigg\downarrow
\Yvcentermath1
\young(1\hfill\hfill\hfill,2\hfill\hfill\hfill)$
\end{center}
\begin{center}
\unitlength 13pt
\begin{picture}(26,5)
\Yboxdim13pt
\put(1,4){\young(\hfill\hfill\hfill\hfill)}
\put(7,1){\young(\hfill\hfill\hfill)}
\put(6.3,1.1){1}
\put(10.2,1.1){1}
\put(7,4){\young(\hfill\hfill\times)}
\put(13,0){\young(\hfill\hfill\hfill,\hfill)}
\put(12.3,0.1){0}
\put(12.3,1.1){1}
\put(14.2,0.1){0}
\put(16.2,1.1){0}
\put(13,4){\young(\hfill\hfill)}
\put(19,0){\young(\hfill\hfill,\hfill)}
\put(18.3,0.1){0}
\put(18.3,1.1){0}
\put(20.2,0.1){0}
\put(21.2,1.1){0}
\put(24,1){\young(\hfill)}
\put(23.3,1.1){0}
\put(25.2,1.1){0}
\end{picture}
\end{center}
\begin{center}
$\Bigg\downarrow
\Yvcentermath1
\young(1\hfill\hfill\hfill,22\hfill\hfill)$
\end{center}
\begin{center}
\unitlength 13pt
\begin{picture}(26,5)
\Yboxdim13pt
\put(1,3){\young(\hfill\hfill\hfill\hfill,\times)}
\put(7,1){\young(\hfill\hfill\hfill)}
\put(6.3,1.1){2}
\put(10.2,1.1){1}
\put(7,4){\young(\hfill\hfill)}
\put(13,0){\young(\hfill\hfill\hfill,\hfill)}
\put(12.3,0.1){0}
\put(12.3,1.1){0}
\put(14.2,0.1){0}
\put(16.2,1.1){0}
\put(13,4){\young(\hfill\hfill)}
\put(19,0){\young(\hfill\hfill,\hfill)}
\put(18.3,0.1){0}
\put(18.3,1.1){0}
\put(20.2,0.1){0}
\put(21.2,1.1){0}
\put(24,1){\young(\hfill)}
\put(23.3,1.1){0}
\put(25.2,1.1){0}
\end{picture}
\end{center}
\begin{center}
$\Bigg\downarrow
\Yvcentermath1
\young(11\hfill\hfill,22\hfill\hfill)$
\end{center}
\begin{center}
\unitlength 13pt
\begin{picture}(26,5)
\Yboxdim13pt
\put(1,4){\young(\hfill\hfill\hfill\hfill)}
\put(7,1){\young(\hfill\hfill\hfill)}
\put(6.3,1.1){1}
\put(10.2,1.1){1}
\put(7,4){\young(\hfill\times)}
\put(13,0){\young(\hfill\hfill\times,\hfill)}
\put(12.3,0.1){0}
\put(12.3,1.1){0}
\put(14.2,0.1){0}
\put(16.2,1.1){0}
\put(13,4){\young(\hfill\hfill)}
\put(19,0){\young(\hfill\hfill,\hfill)}
\put(18.3,0.1){0}
\put(18.3,1.1){0}
\put(20.2,0.1){0}
\put(21.2,1.1){0}
\put(24,1){\young(\hfill)}
\put(23.3,1.1){0}
\put(25.2,1.1){0}
\end{picture}
\end{center}
\begin{center}
$\Bigg\downarrow
\Yvcentermath1
\young(11\hfill\hfill,223\hfill)$
\end{center}
\begin{center}
\unitlength 13pt
\begin{picture}(26,5)
\Yboxdim13pt
\put(1,3){\young(\hfill\hfill\hfill\hfill,\times)}
\put(7,1){\young(\hfill\hfill\times)}
\put(6.3,1.1){1}
\put(10.2,1.1){1}
\put(7,4){\young(\hfill)}
\put(13,0){\young(\hfill\hfill,\hfill)}
\put(12.3,0.1){0}
\put(12.3,1.1){0}
\put(14.2,0.1){0}
\put(15.2,1.1){0}
\put(13,4){\young(\hfill\hfill)}
\put(19,0){\young(\hfill\hfill,\hfill)}
\put(18.3,0.1){0}
\put(18.3,1.1){0}
\put(20.2,0.1){0}
\put(21.2,1.1){0}
\put(24,1){\young(\hfill)}
\put(23.3,1.1){0}
\put(25.2,1.1){0}
\end{picture}
\end{center}
\begin{center}
$\Bigg\downarrow
\Yvcentermath1
\young(112\hfill,223\hfill)$
\end{center}
\begin{center}
\unitlength 13pt
\begin{picture}(26,5)
\Yboxdim13pt
\put(1,4){\young(\hfill\hfill\hfill\hfill)}
\put(7,1){\young(\hfill\hfill)}
\put(6.3,1.1){1}
\put(9.2,1.1){1}
\put(7,4){\young(\times)}
\put(13,0){\young(\hfill\hfill,\times)}
\put(12.3,0.1){0}
\put(12.3,1.1){0}
\put(14.2,0.1){0}
\put(15.2,1.1){0}
\put(13,4){\young(\hfill\hfill)}
\put(19,0){\young(\hfill\hfill,\times)}
\put(18.3,0.1){0}
\put(18.3,1.1){0}
\put(20.2,0.1){0}
\put(21.2,1.1){0}
\put(24,1){\young(\times)}
\put(23.3,1.1){0}
\put(25.2,1.1){0}
\end{picture}
\end{center}
\begin{center}
$\Bigg\downarrow
\Yvcentermath1
\young(112\hfill,2235)$
\end{center}
\begin{center}
\unitlength 13pt
\begin{picture}(26,5)
\Yboxdim13pt
\put(1,2){\young(\hfill\hfill\hfill\hfill,\times)}
\put(7,0){\young(\hfill\times)}
\put(6.3,0.1){1}
\put(9.2,0.1){1}
\put(7.5,3){$\emptyset$}
\put(13,0){\young(\hfill\times)}
\put(12.3,0.1){0}
\put(15.2,0.1){0}
\put(13,3){\young(\hfill\hfill)}
\put(19,0){\young(\hfill\times)}
\put(18.3,0.1){0}
\put(21.2,0.1){0}
\put(24,0){$\emptyset$}
\end{picture}
\end{center}
\begin{center}
$\Bigg\downarrow
\Yvcentermath1
\young(1124,2235)$
\end{center}
\begin{center}
\unitlength 13pt
\begin{picture}(26,5)
\Yboxdim13pt
\put(1,3){\young(\hfill\hfill\hfill\hfill)}
\put(7,0){\young(\hfill)}
\put(6.3,0.1){0}
\put(8.2,0.1){0}
\put(7.5,3){$\emptyset$}
\put(13,0){\young(\hfill)}
\put(12.3,0.1){0}
\put(14.2,0.1){0}
\put(13,3){\young(\hfill\hfill)}
\put(19,0){\young(\hfill)}
\put(18.3,0.1){0}
\put(20.2,0.1){0}
\put(24,0){$\emptyset$}
\end{picture}
\end{center}
Continuing in this way, we obtain the following path:
$$
\Yvcentermath1
\young(1111)\otimes
\young(12,23,34)\otimes
\young(1124,2235)\, .$$
In order to make comparison, we consider
the same example in Example \ref{ex:led}
in the main text.

\end{example}

\subsection{Several facts about KSS bijection}
\label{app:facts}
As we have seen before,
the map $\phi^{-1}$ is a function of
choice of rows of the quantum space.
Its dependence is described by the following
theorem.

\begin{theorem}[\cite{KSS}, Lemma 8.5]\label{th:KSS}
Take two rows from the quantum space
of the rigged configuration arbitrarily,
and denote them by $\nu_a$ and $\nu_b$.
When we remove successively these two rows as
$\nu_a$ at first
and next $\nu_b$ under $\phi^{-1}$,
then we obtain two tableaux,
which we denote by $a_1$ and $b_1$, respectively.
Next, on the contrary, we first remove $\nu_b$
and second $\nu_a$
$($keeping the order of other removal invariant$)$
and we get $b_2$ and $a_2$.
Then we have
\begin{equation}
b_1\otimes a_1\,\simeq\, a_2\otimes b_2 ,
\end{equation}
under the isomorphism of the
combinatorial $R$.
\end{theorem}

Recently, the original KSS bijection was
extended to include non-highest weight elements
\cite{Sch,DS}.
The combinatorial procedures for $\phi$
and $\phi^{-1}$ are the formal extension
of the original ones.
For arbitrary element
$b\in B^{r_1,s_1}\otimes
B^{r_2,s_2}\otimes\cdots\otimes B^{r_L,s_L}$,
we obtain $\phi (b)$ as an extension of
the rigged configuration, which we call
{\it unrestricted rigged configurations}.
Characterization for the set of the
unrestricted rigged configurations is
given in these papers, however
we omit it since we do not use it.
For our purpose, it is enough to
start from arbitrarily given path,
and use the prescription given in
Section 7.1 of \cite{KSY}.
In the prescription given there,
we put suitable prefix to the path,
and the effect of such modification
of path on the level of the rigged
configuration is explicitly estimated.
As the result, many apparatus
(like Theorem \ref{th:KSS}) holds
for the case of the unrestricted rigged
configuration.
In particular, our procedure given in
Theorem \ref{th:main} does not depend
on whether the path is highest
or non-highest.

\bigskip

\noindent
{\bf Acknowledgements:}
The author is grateful to Atsuo Kuniba
for collaboration at an
early stage of the present study,
and to Yasuhiko Yamada for careful reading
of the manuscript.
He is a research fellow of the 
Japan Society for the Promotion of Science.

\end{document}